\documentclass{article}

\usepackage[top=1 in, bottom=1 in, left=1 in, right = 1 in]{geometry}
\usepackage{amsmath,enumerate}
\usepackage{amsthm}
\usepackage{amssymb}
\usepackage{setspace}
\usepackage{graphicx,float,subfloat,epstopdf,subfig}
\usepackage[colorlinks=true,urlcolor=blue,citecolor=OliveGreen,linkcolor=red]{hyperref}
\graphicspath{ {images_SIRS/} }
\usepackage[dvipsnames]{xcolor}
\usepackage{algorithm}
\usepackage{algpseudocode}
\usepackage[font=small,labelfont=bf,width=0.9\textwidth]{caption}

\newcounter{thmCounter}[subsection]

\newcommand{\dref}[1]{definition~\hyperref[def:#1]{\color{red}#1}}
\newcommand{\tref}[1]{theorem~\hyperref[thm:#1]{\color{red}#1}}
\newcommand{\Dref}[1]{Definition~\hyperref[def:#1]{\color{red}#1}}
\newcommand{\Tref}[1]{Theorem~\hyperref[thm:#1]{\color{red}#1}}
\newcommand{\lref}[1]{lemma~\hyperref[lem:#1]{\color{red}#1}}
\newcommand{\Lref}[1]{Lemma~\hyperref[lem:#1]{\color{red}#1}}
\newcommand{\pref}[1]{proposition~\hyperref[prop:#1]{\color{red}#1}}
\newcommand{\Pref}[1]{Proposition~\hyperref[prop:#1]{\color{red}#1}}
\newcommand{\eref}[1]{example~\hyperref[ex:#1]{\color{red}#1}}
\newcommand{\Eref}[1]{Example~\hyperref[ex:#1]{\color{red}#1}}
\newcommand{\corref}[1]{corollary~\hyperref[ex:#1]{\color{red}#1}}
\newcommand{\Corref}[1]{Corollary~\hyperref[ex:#1]{\color{red}#1}}

\newcounter{count}
\stepcounter{count}

\DeclareMathOperator*{\argmin}{arg\,min}

\newcommand{\abs}[1]{\left| #1  \right|}

\newcommand{\R}{\mathbb{R}}

\newcommand{\W}{\mathcal W}

\newcommand{\QED}{\hfill $\blacksquare$ \\ }

\newcommand{\defn}{\refstepcounter{thmCounter}\noindent {\bf Definition \arabic{section}.\arabic{thmCounter}. }\label{def:\arabic{section}.\arabic{thmCounter}}}
\newcommand{\prf}{\noindent {\bf Proof. }}
\newcommand{\cor}{\refstepcounter{thmCounter} \noindent {\bf Corollary \arabic{section}.\arabic{thmCounter}. }\label{thm:\arabic{section}.\arabic{thmCounter}}}
\newcommand{\thm}{\refstepcounter{thmCounter} \noindent {\bf Theorem \arabic{section}.\arabic{thmCounter}. }\label{thm:\arabic{section}.\arabic{thmCounter}}}
\newcommand{\lem}{\refstepcounter{thmCounter} \noindent {\bf Lemma \arabic{section}.\arabic{thmCounter}. }\label{lem:\arabic{section}.\arabic{thmCounter}}}

\newcommand{\weakto}{\rightharpoonup}

\newcommand*{\defeq}{\mathrel{\vcenter{\baselineskip0.5ex \lineskiplimit0pt
                     \hbox{\scriptsize.}\hbox{\scriptsize.}}}%
     		     =}

\def\XXint#1#2#3{{\setbox0=\hbox{$#1{#2#3}{\int}$ }
\vcenter{\hbox{$#2#3$ }}\kern-.6\wd0}}

\newtheorem{theorem}{Theorem}[section]

\title{Optimal Feedback Control of the Fokker-Planck Equation for an SIRS Epidemic Model}
\author{Christian Parkinson, Souvik Roy}
\date{}
\begin{document}

\maketitle

\begin{abstract}
We develop a Fokker--Planck feedback control framework for a
stochastic susceptible--infected--recovered--susceptible (SIRS)
epidemic model with three intervention mechanisms targeting
disease transmission, loss of immunity, and recovery. The
interventions are formulated as space- and time-dependent
feedback controls that influence both the drift and diffusion
coefficients of the underlying stochastic differential equations,
allowing control of the probability distribution of epidemic
states. We investigate the existence of optimal controls and
derive first-order necessary optimality conditions involving
a forward Fokker--Planck equation, a backward adjoint equation,
and a Hamiltonian minimization condition. For numerical
approximation, we combine a Chang--Cooper discretization of
the diagonal drift--diffusion terms with a finite-difference
treatment of cross-diffusion and a sequential quadratic
Hamiltonian (SQH) algorithm for computing the optimal feedback
policies. We also examine sufficient-decrease and convergence
properties of the SQH iteration. Numerical experiments consider
four epidemiological objectives involving cumulative infection
burden, infection-threshold exceedance, terminal susceptible
population, and a combination of running and terminal costs.
The resulting feedback policies exhibit distinct spatial and
temporal structures depending on the epidemiological objective.
Monte Carlo simulations of the underlying stochastic system
illustrate that the optimized interventions can substantially
reduce infection growth and undesirable epidemic outcomes
relative to uncontrolled dynamics. The proposed framework
provides a distribution-based approach to stochastic epidemic
control that integrates mathematical analysis, numerical
optimization, and state-dependent intervention design.
\end{abstract}



\section{Introduction}\label{sec:intro}

Mathematical modeling has long played an important role in
understanding the transmission of infectious diseases and
evaluating strategies for epidemic prevention and control.
Beginning with the classical work of Kermack and McKendrick
\cite{Kermack1927}, compartmental epidemic models have
provided a systematic framework for describing the evolution
of susceptible, infected, and recovered populations. Among
these, the susceptible--infected--recovered--susceptible
(SIRS) model is particularly relevant to infectious diseases
for which immunity is temporary or incomplete. By allowing
recovered individuals to return to the susceptible class,
the SIRS model accounts for reinfection and the possibility
of recurrent epidemic outbreaks. Such dynamics motivate
the development of intervention strategies that reduce
transmission, limit the loss of immunity, and enhance
recovery while accounting for the costs associated with
implementing these measures.

Although deterministic epidemic models provide valuable
insights into disease transmission, their predictions are
based on prescribed parameter values and initial conditions.
In practice, epidemic dynamics are influenced by substantial
uncertainty arising from variations in contact patterns,
transmission rates, recovery mechanisms, immunity duration,
and environmental or behavioral factors. Stochastic epidemic
models address these uncertainties by incorporating random
perturbations into the governing dynamics. Several authors investigated stochastic epidemic
dynamics using stochastic differential equations and
established important results concerning positivity,
extinction, persistence, stationary distributions and the effects of different types of noise \cite{tornatore,Gray,lanconelli,PW3,tornatore2}. See \cite{allen2017primer} for an overview of such models. 
Generally, stochastic differential equation models
describe ensembles of possible epidemic trajectories
rather than a single deterministic evolution. This
distinction is particularly important when intervention
strategies must account for the probability of severe
outbreaks or other undesirable epidemiological outcomes.

Optimal control theory provides a mathematical framework
for determining intervention strategies that balance
epidemiological benefits against economic, logistical,
and societal costs. In compartmental epidemic models,
controls may represent nonpharmaceutical interventions
that reduce transmission, vaccination or other measures
that influence immunity, and treatment strategies that
accelerate recovery \cite{OC4,OC1}, and researchers have analyzed models incorporating confounding effects ranging from noncompliance with control measures \cite{BPW,PW4} to media coverage \cite{OC3} to limited resources \cite{OC5}. Classical optimal control formulations
typically determine time-dependent interventions by
minimizing a prescribed objective functional subject to
the epidemic dynamics. In stochastic settings, optimal
control can also be formulated through controlled diffusion
processes and dynamic programming, as developed in the
general theory of stochastic control
\cite{FlemingSoner2006}. However, time-dependent
interventions do not explicitly adapt to the current
epidemic state. Moreover, when the objective involves
the likelihood of undesirable outcomes rather than
only expected population levels, it is useful to
consider the evolution of the entire probability
distribution of epidemic states.

The Fokker--Planck (FP) equation offers a natural framework
for this purpose. Associated with an underlying stochastic
differential equation, the FP equation describes the temporal evolution
of the probability density of the stochastic state
variables. Consequently, an optimal control problem
constrained by the FP equation allows
interventions to influence the distribution of possible
epidemic outcomes. This density-based perspective
provides a direct way to incorporate objectives involving
expected infection burden, probability mass in
high-infection regions, and desirable terminal epidemic
states. Furthermore, by allowing the controls to depend
on both the epidemic state and time, the resulting
formulation accommodates feedback intervention policies.
Unlike purely time-dependent controls, these policies
can assign different intervention intensities to
different regions of the epidemic state space. While FP control has been applied biological and ecological models such as those for pedestrian motion \cite{RBH} and biological imaging \cite{AB}, it has not been broadly applied to epidemiological models with the exception of \cite{PS} where time-dependent controls were considered. Indeed, most work regarding FP control assumes time-dependent control or specially parameterized control ansatzes and, notably, does not consider control of the diffusive terms in the FP equation \cite{BB,OptControlFPSpatial}.  

Despite the modeling advantages, FP optimal control
problems introduce mathematical and computational
challenges. The probability density satisfies a
time-dependent partial differential equation whose
coefficients depend on the control variables. When
the controls influence both the drift and diffusion
coefficients, the resulting control-to-state mapping
is highly nonlinear, and the diffusion operator may contain
mixed spatial derivatives. Establishing existence of
optimal controls and deriving necessary optimality
conditions therefore requires an appropriate analytical
framework for the controlled FP equation.
The numerical approximation must also account for
the nonnegativity and conservation properties of
probability densities. The classical Chang--Cooper
method \cite{ChangCooper1970} provides a useful
structure-preserving discretization for drift--diffusion
equations, although additional treatment is required
for multidimensional systems with cross-diffusion.

Motivated by these considerations, we develop a
FP optimal control framework for a stochastic
SIRS epidemic model with three distinct intervention
mechanisms. Stochastic perturbations are introduced
into transmission, loss of immunity, and recovery,
and the corresponding feedback controls modify both
the deterministic transition rates and the stochastic
fluctuations. The first control reduces transmission
through nonpharmaceutical interventions, the second
reduces loss of immunity, and the third increases
the recovery rate. The resulting stochastic system
induces a two-dimensional FP equation
with a control-dependent diffusion matrix and
nonzero cross-diffusion terms.

The principal contributions of this work are as follows.
First, we formulate a distribution-based optimal control
problem with space- and time-dependent feedback controls
and a general objective functional incorporating
intervention costs, running epidemic costs, and terminal
costs. In particular, our problem includes controls in the diffusion terms of the FP equation. Second, we investigate the analytical properties
of the associated control-to-state operator and develop
existence and first-order optimality results using a
forward FP equation, a backward adjoint
equation, and a Hamiltonian minimization condition.
Third, we construct a numerical framework combining
Chang--Cooper discretization of the diagonal
drift--diffusion terms, finite-difference approximation
of the mixed diffusion term, and a sequential quadratic
Hamiltonian (SQH) algorithm for computing feedback
controls. We also investigate sufficient-decrease
and convergence properties of the SQH iteration. 

The numerical experiments consider four epidemiological
objectives that demonstrate the flexibility of the
proposed framework. The first minimizes cumulative
infection burden, while the second penalizes epidemic
states in which the infected fraction exceeds a
prescribed threshold. The third introduces a terminal
reward for maintaining a sufficiently large susceptible
population, and the fourth combines threshold-based
infection regulation with a terminal susceptible-state
objective. For each scenario, we compute the evolving
probability density and the corresponding feedback
control maps. Monte Carlo simulations of the underlying
controlled and uncontrolled stochastic differential
equations are then used to examine the effects of
the optimized interventions on individual epidemic
trajectories and their ensemble behavior.

The remainder of the paper is organized as follows.
Section~\ref{sec:model} introduces the controlled
stochastic SIRS model, derives the associated
FP equation, and formulates the feedback
optimal control problem. Section~\ref{sec:analysis}
develops the analytical framework, including existence
of an optimal control and first-order necessary
optimality conditions. Section~\ref{sec:computation}
presents the numerical approximation of the forward
FP and backward adjoint equations together
with the SQH algorithm and its convergence analysis.
Section~\ref{sec:results} presents numerical experiments
under different epidemiological objectives and compares
the resulting controlled and uncontrolled stochastic
dynamics.

\section{A controlled stochastic SIRS model and FP equation}
\label{sec:model}

We begin from the following controlled Susceptible-Infected-Recovered-Susceptible (SIRS) dynamics which model the spread of a generic disease in a generic region: \begin{equation} \label{eq:SIRSdeterministic}\begin{aligned} S' &= b - (1-u_1) \beta SI + (1-u_2)\lambda R -\delta S, \\
I' &= (1-u_1)\beta SI - ((1+u_3)\gamma+\delta) I, \\
R' &= (1+u_3)\gamma I -((1-u_2)\lambda + \delta)R.
\end{aligned}\end{equation} Here $b$ represents the constant introduction of new members into the population, and $\delta$ represents the natural death rate. A basic assumption of this model is that the population is homogeneous and well-mixed so that infection spreads via mass action with rate $\beta$ and infected individuals recover at rate $\gamma$. The parameter $\lambda$ represents loss of immunity whereby recovered individuals become susceptible again. We supply the policy-maker with three control variables $u \defeq (u_1,u_2,u_3)$ by which they may inhibit disease spread: \begin{enumerate}[(a)]
\item $u_1 \in [0,u_1^{\text{max}}]$ for some $u_1^{\text{max}} \in [0,1]$ - this is reduction in infectivity due to governmentally enacted non-pharmaceutical intervention measures,
\item $u_2 \in [0,u_2^{\text{max}}]$ for some $u_2^{\text{max}} \in [0,1]$ - reduction in loss of immunity (perhaps due to vaccination or other treatment efforts),
\item $u_3 \in [0,u_3^{\text{max}}]$ for some $u_3^{\text{max}} \in [0,1]$ - increase in recovery rate due to treatment efforts.
\end{enumerate} Due to the highly volatile nature of epidemic disease spread, we also add stochasticity to the model, and specifically following \cite{Gray,PW3}, we assume that there is some uncertainty in the rate parameters $\beta,\gamma,\lambda$. Thus we formally replace $\beta$ with $\beta + \sigma_\beta \dot W$, and likewise with $\gamma,\lambda$. Finally, we assume that $b = \delta$ and that the initial conditions are normalized so that $S(0)+I(0)+R(0)=1$. Under these conditions, we see that the total population $N = S+I+R$ described by \eqref{eq:SIRSdeterministic} (and indeed by the stochastic generalization) will be constant: $N(t) \equiv 1$. Thus we can write $R = 1-S-I$, and drop this compartment. With all these assumptions, we introduce the controlled, stochastic SIRS model:\begin{equation}\label{eq:dynamicsControlled}
\begin{aligned} dS &= (b - (1-u_1)\beta SI + (1-u_2)\lambda(1-S-I)- \delta S)dt\\ &\hspace{4cm} -(1-u_1)\sigma_\beta SIdW_1 + (1-u_2)\sigma_\lambda (1-S-I)dW_2, \\ 
dI &= ((1-u_1)\beta SI - (1+u_3)\gamma I - \delta I)dt + (1-u_1)\sigma_\beta SI dW_1  -(1+u_3)\sigma_\gamma IdW_3.
\end{aligned} \end{equation} We write the dynamics in the general form in the form \begin{equation} \label{eq:dynamicsGeneral} dx = F(x,u)dt + \Sigma(x,u) dW \end{equation} where $x = (S,I)$, $u = (u_1,u_2,u_3)$, and $W$ is a 3 dimensional Brownian motion with independent components. In this notation, the drift and diffusion functions are \begin{align*} F(x,u) =& \hphantom{\frac 1 2}\begin{bmatrix} b - (1-u_1)\beta SI + (1-u_2)\lambda(1-S-I)- \delta S \\ (1-u_1)\beta SI - (1+u_3)\gamma I - \delta I\end{bmatrix}, \\
\frac1 2 \Sigma(x,u)\Sigma(x,u)^T =& \frac 1 2\begin{bmatrix}
    ((1-u_1)\sigma_\beta SI)^2 + ((1-u_2)\sigma_\lambda (1-S-I))^2 & -((1-u_1)\sigma_\beta SI)^2 \\ -((1-u_1)\sigma_\beta SI)^2 & ((1-u_1)\sigma_\beta SI)^2 + ((1+u_3)\sigma_\gamma I)^2
\end{bmatrix}. \end{align*}

For notational convenience, we set $D(x,u) = \tfrac 1 2 \Sigma(x,u) \Sigma(x,u)^T$. We record some important observations regarding the model. 
\begin{itemize}
    \item By the arguments in \cite{Gray,PW3}, if $S_0,I_0,R_0 > 0$ and $S_0+I_0+R_0 \le 1$, then $[0,1]^3$ is a forward invariant set for the dynamics, meaning that we can take $\Omega = [0,1]^2$ and let $x_0 = (S_0,I_0) \sim \rho_0$ where $\rho_0$ is supported in $\Omega$. 
    \item Our admissible control actions take values in the rectangle $R_U = \prod_{i=1}^3 [0, u^{\text{max}}_i]$.
    \item All the components of the functions $F:\Omega \times R_U \to \R^2$ and $D:\Omega\times R_U \to \R^{2\times 2}$ and their $x$-derivatives are continuous in $(x,u)$ and thus bounded on $\Omega \times R_U$.
\end{itemize}

The Fokker-Planck equation corresponding to \eqref{eq:dynamicsGeneral} is \begin{equation}\label{eq:controlledFP} \begin{aligned}
&\rho_t + \nabla \cdot(F(x,u)\rho) = \sum_{i,j=1}^2 \frac{\partial^2}{\partial x_i \partial x_j} (D_{ij}(x,u)\rho), \,\,\,\,\, \text{ in } \Omega \times (0,T], \\
&\rho(x,0) = \rho_0(x), \,\,\,\, \text{ in } \,\, \Omega \\
&E(x,u)\cdot n = 0 \,\, \text{ on } \,\, \partial \Omega \times (0,T],
\end{aligned}
\end{equation} where the flux $E(x,u)$ is defined componentwise by \begin{equation} \label{eq:FLUX}E_i(x,u) = F_i(x,u)\rho - \sum^2_{j=1} \frac{\partial}{\partial x_j}(D_{ij}(x,u)\rho), \,\,\,\,\,\,\,\,\, i =1,2.\end{equation}

We now formulate an optimal control problem. We are interested in feedback control, so we consider $u_i : \Omega_T\to R_U$ for each $i=1,2,3$ where $\Omega_T = \Omega \times[0,T]$. Note that because the controls enter into \eqref{eq:controlledFP} in the diffusion coefficients, we need very regular controls to guarantee nice properties of the solution. Some of this additional regularity will be required for the ensuing proof of existence of minimizers for the optimal control problem in \tref{3.5} (and the preceding lemmas), and yet more is required for the proof of the Pontryagin Minimum Principle in \tref{3.8}.  Accordingly, we take $$U_{\text{ad}} = \{u \in W^{1,\infty}(\Omega_T;R_U) \,\, | \,\,  \|u_i\|_{W^{1,\infty}(\Omega_T)} \le M, \,\,\,\, i = 1,2,3\}.$$ for a fixed positive constant $M$. We consider a relatively generic cost functional \begin{equation} \label{eq:costFunc}J(\rho,u) = \int_{0}^T \int_\Omega \Big( \ell(u(x,t)) + G(x,t)\rho(x,t)\Big) dxdt + \int_\Omega K(x)\rho(x,T)dx\end{equation} where $\ell:R_U \to [0,\infty)$ and $G,K$ are, at a minimum, bounded. For particular theorems, these functions will be required to satisfy further conditions as specified below.  

Our optimization problem can then be stated abstractly as follows \begin{equation}\label{eq:optProb}
\inf_{u \in U_{\text{ad}}} \Big\{ J(\rho,u) \,\, : \,\, \rho \text{ satisfies \eqref{eq:controlledFP} with control } u\Big\} 
\end{equation}

Here the functions $\ell, G, K$ express the preferences of the policy-maker. The control cost $\ell(u)$ represents represents any political, economic, of logistical costs associated with implementing control measures. For example, this could encode the loss of gross domestic product (GDP) incurred when requiring businesses to close or the costs associated with manufacturing and distributing vaccines. The running cost $G(x,t)$ and final cost $K(x)$ implicitly define the desirable or undesirable epidemic outcomes. For example, recalling that $x = (S,I)$, a policy-maker may set $G(x,t) = I$, meaning that more cost is assessed when there is more probability density associated with states where the infectious population is larger.  

\section{Analysis of the FP optimal control problem}
\label{sec:analysis}
 We prove existence of a minimizer for \eqref{eq:optProb} via analysis of the control-to-state map. Before this, we need a basic existence result. \\

\defn For a fixed control function $u \in U_{\text{ad}}$ and initial condition $\rho_0 \in L^2(\Omega)$, a function $\rho(x,t)$ is formally defined to be a weak solution of \eqref{eq:controlledFP} if it satisfies $\rho(\cdot,0) = \rho_0$ and \begin{equation} \label{eq:weakSoln}\int_\Omega \bigg(\rho_t \varphi - \rho F(x,u)\cdot \nabla \varphi + \sum^2_{i,j=1} \frac{\partial}{\partial x_j} (D_{ij}(x,u) \rho)\frac{\partial \varphi}{\partial x_i}\bigg)dx = 0\end{equation} for all $\varphi \in H^1(\Omega)$ and a.e. $t \in [0,T]$. \\

\thm Fix $\rho_0 \in H^1(\Omega)$. Then for any $u \in U_{\text{ad}}$ there is a unique weak solution $\rho$ of \eqref{eq:controlledFP} with $$\rho \in L^2(0,T;H^1(\Omega)) \,\,\,\, \text{ and } \,\,\,\, \rho_t \in L^2(0,T; H^{-1}(\Omega)).$$ Moreover, $\|\rho\|_{L^2(0,T;H^1(\Omega))}$ and $\|\rho_t\|_{L^2(0,T;H^{-1}(\Omega))}$ are uniformly bounded over choices of $u \in U_{\text{ad}}.$ Further, $$\int_\Omega \rho(x,t)dx = \int_\Omega \rho_0(x)dx, \,\,\,\,\,\, \text{ for all } \, t \in (0,T]$$ and if $\rho_0 \ge 0$ almost everywhere in $\Omega$, then $\rho(\cdot,t) \ge 0$ almost everywhere in $\Omega$ for all $t \in (0,T].$ \\

The proof of this theorem is classical so we omit it, but comment on some facets. Existence and uniqueness of the solution follows as in \cite[Ch. 7.1]{Evans}. The uniform boundedness over choices of controls follows since $F(x,u)$ and $D(x,u)$ and their spatial derivatives are uniformly bounded for $(x,u) \in \Omega \times R_U$. Conservation of the integral follows by taking $\varphi \equiv 1$ in \eqref{eq:weakSoln}, and nonnegativity follows by setting $\varphi = (\rho(\cdot,t))^- \defeq -\min(\rho(\cdot,t),0)$ at each time $t$ and establishing the bound $\frac d {dt}\|\rho(\cdot,t)^-\|_{L^2(\Omega)} \le C \|\rho(\cdot,t)^-\|_{L^2(\Omega)}$, whereupon the result follows by Gr\"onwall's inequality since $\rho_0^- = 0$. Note that because $\rho_0$ is a probability density function, we will henceforth always assume that $\rho_0 \ge 0$ and $\|\rho_0\|_{L^1(\Omega)} = 1$, whereupon we have these same properties for $\rho(\cdot,t)$ for all times $t \in [0,T]$. We note that \tref{3.2} does not require much regularity for the controls $u$, so for this theorem the regularity in our control space $U_{\text{ad}}$ is sufficient but certainly unnecessary.

From \tref{3.2}, we have a well defined control-to-state operator \begin{equation} \label{eq:CtoS} \begin{aligned}
    \mathcal S : U_{\text{ad}} &\to \W\defeq \{ \rho \in L^2(0,T;H^1(\Omega)) \,\, : \,\, \rho_t \in L^2(0,T; H^{-1}(\Omega))\} \subset L^2(\Omega_T),\\
    u \in U_{\text{ad}} &\mapsto \rho \in \W, \text{ the weak solution of } \eqref{eq:controlledFP} \text{ with control map } u.
    \end{aligned}
\end{equation} We prove existence of a solution of \eqref{eq:optProb} through analysis of $\mathcal S$. We proceed with two lemmas and then state existence of minimizers as a theorem. In all of the following work, we equip $U_{\text{ad}}$ with the $L^2$-topology induced by the norm $$\|u\|_{L^2(\Omega_T)} = \left(\sum^3_{i=1} \|u_i\|_{L^2(\Omega_T)}^2\right)^{1/2}.$$ Additionally, whenever a subsequence with some property is mentioned, we implicitly pass to this subsequence, but to avoid excessive notation, we do not rename the sequences.\\

\lem The control-to-state operator $\mathcal S$ is weakly closed. \\

\prf Take a sequence of controls $u_k \in U_{\text{ad}}$ such that $u_k \weakto u$ in $L^2(\Omega_T)$ and $\rho_k \defeq \mathcal S(u_k) \weakto \rho$ in $\mathcal W$. We need to prove that $u \in U_{\text{ad}}$ and $\mathcal S(u) = \rho.$ 

By the Banach-Alaoglu theorem and compatibility of weak limits, along a subsequence, we can achieve $u_k \rightharpoonup^* u \in L^\infty(\Omega_T;R_U)$ and $(\nabla u_i)_k \rightharpoonup^* g_i$ in $L^\infty(\Omega_T;\R^2)$ for some $g_i \in L^\infty(\Omega_T;\R^2)$, $i = 1,2,3$. Now for any $\varphi \in C^\infty_c(\Omega_T)$, $i=1,2,3$ and any first order partial derivative $\partial$, we see $$\int_{\Omega_T} (u_i)_k \partial \varphi dx dt = -\int_{\Omega_T} \partial (u_i)_k \varphi dxdt.$$ By weak$^*$ convergence of $u_k$ to $u$ and of $(\nabla u_i)_k$ to $g_i$, this gives $$\int_{\Omega_T} u_i \partial \varphi dx dt = - \int_{\Omega_T} \partial(g_i) \varphi dx dt, \,\,\,\,\,\,\, \text{ for all } \,\, \varphi \in C^\infty_c(\Omega_T)$$ implying that $u_i$ is weakly differentiable with weak partial $\partial (g_i)$, and $$\|\partial u_i\|_{L^\infty(\Omega_T)} \le \liminf_{k\to\infty} \|\partial (u_i)_k\|_{L^\infty(\Omega_T)} \le M$$ by weak$^*$ lower semicontinuity of the $L^\infty$ norm. Finally, weak $L^2$ convergence also preserves pointwise a.e. bounds, as follows by testing $u_k$ against appropriately designed indicator functions, and thus $u(x,t) \in R_U$ for almost every $(x,t) \in \Omega_T$. We conclude that $u \in U_{\text{ad}}.$ 

It remains to prove that $\mathcal S(u) = \rho$. By \tref{3.2}, the sequence $\rho_k$ is bounded in $\W$. Since we have compact-continuous embeddings $$H^1(\Omega) \hookrightarrow\hookrightarrow L^2(\Omega) \hookrightarrow H^{-1}(\Omega),$$ the Aubin-Lions lemma \cite[Ch. III, Prop. 1.3]{Showalter} yields a subsequence of $\rho_k$ which converges strongly to $\rho$ in $L^2(\Omega_T).$ Since $u_k(x,t) \in R_U$ for a.e. $(x,t) \in \Omega_T$ and $\|u_k\|_{W^{1,\infty}(\Omega_T)} \le M$, the sequence $u_k$ is uniformly bounded and uniformly Lipschitz, whereupon the Arzela-Ascoli theorem gives $u_k \to u$ uniformly along a subsequence. Because $F$ is continuous and uniformly bounded, this allows us to conclude that $F(x,u_k) \to F(x,u)$ in $L^2(\Omega_T)$ via the dominated convergence theorem (and similarly for $D_{ij}$ and their spatial derivatives). Finally to prove that $\mathcal S(u) = \rho$, we fix $\varphi \in H^1(\Omega)$. Considering \eqref{eq:weakSoln}, there are three terms to address. The pairing $(\rho_k)_t \varphi$ behaves nicely since $\rho_k\weakto \rho$ in $\W$. The term $\rho_kF(x,u_k)\cdot \nabla \varphi$ is controlled since both $\rho_k \to \rho$ and $F(x,u_k) \to F(x,u)$ strongly in $L^2(\Omega_T)$. And upon expanding, the terms $\frac{\partial}{\partial x_j} (D_{ij}(x,u_k) \rho_k)\frac{\partial \varphi}{\partial x_i}$ are controlled since $D_{ij}(x,u_k)$ is uniformly bounded and converges strongly to $D_{ij}(x,u)$ in $L^2(\Omega_T)$, and $\rho_k \to \rho$ weakly in $\W$. Thus we find that $\rho$ is a weak solution of \eqref{eq:controlledFP} corresponding to $u$; that is, $\mathcal S(u) = \rho$. \QED

\lem The control-to-state operator $\mathcal S$ is weak-to-strong continuous in the $L^2$-topology. That is, if $u_k \in U_{\text{ad}}$ and $u_k \weakto u$ in $L^2(\Omega_T)$, then $\mathcal S(u_k) \to S(u)$ in $L^2(\Omega_T)$.\\

\prf Suppose $u_k \in U_{\text{ad}}$ and $u_k \weakto u$ in $L^2(\Omega_T)$. Once again, $\rho_k \defeq \mathcal S(u_k)$ is a bounded sequence in $\W$, and thus any subsequence is bounded in $\W$. Hence by the Aubin-Lions lemma \cite[Ch. III, Prop. 1.3]{Showalter}, we can pass to a further subsequence which converges strongly to some $\rho \in L^2(\Omega_T)$. By the closure guaranteed in \lref{3.3}, we must have $\rho = \mathcal S(u)$. Hence we have prove that every subsequence of $\mathcal S(u_k)$ has a further subsequence which converges strongly to $\mathcal S(u)$ in $L^2(\Omega_T)$. By the Urysohn subsequence principle \cite[\S2.1.17]{TaoMeasure}, we have that $\mathcal S(u_k) \to \mathcal S(u)$ in $L^2(\Omega_T)$. \QED

We recall before stating our existence theorem that the functions $\ell, G, K$ are used in the definition of the cost functional $J(\rho,u)$ in \eqref{eq:costFunc}.\\

\thm Assume that $\ell:R_U \to [0,\infty)$ is convex and lower-semicontinuous, $G \in L^\infty(\Omega_T)$, and $K \in L^\infty(\Omega)$. Then there exists a solution $u \in U_{\text{ad}}$---and corresponding weak solution $\rho \in \mathcal W$ of \eqref{eq:controlledFP}---of the minimization problem \eqref{eq:optProb}. \\

\prf Define $$m \defeq \inf_{u \in U_{\text{ad}}} \Big\{ J(\rho,u) \,\, : \,\, \rho \text{ satisfies \eqref{eq:controlledFP} with control } u\Big\}.$$ Recall that $\rho(\cdot,t) \ge 0$ and $\|\rho(\cdot,t)\|_{L^1(\Omega)} = 1$ for all $t \in [0,T]$. Then whenever $\ell,G,K$ satisfy the given assumptions, $$m \ge -T\abs{\Omega}\|G\|_{L^\infty(\Omega_T)} - \abs{\Omega}\|K\|_{L^\infty(\Omega)} > -\infty.$$ Thus we can take a sequence $u_k \in U_{\text{ad}}$ with corresponding $\rho_k \defeq \mathcal S(u_k) \in W$ such that $J(\rho_k,u_k) \searrow m$. Since $u_k$ is, in particular, a bounded sequence in $L^2(\Omega_T)$, there exists a subsequence converging weakly to some $u \in L^2(\Omega_T)$ and by the weak closure in \lref{3.3}, we have $u \in U_{\text{ad}}.$ By the weak-to-strong continuity in \lref{3.4}, we have $\rho_k \to \rho \defeq \mathcal S(u)$ in $L^2(\Omega_T)$. Considering the definition of the cost functional in \eqref{eq:costFunc}, the conditions that $\ell$ is nonnegative, convex and lower-semicontinuous guarantee that the map $u \in U_{\text{ad}} \mapsto \int_{\Omega_T} \ell(u) dxdt \in [0,\infty)$ is weakly lower-semicontinuous, and the integral terms involving the $\rho_k$ are preserved under limits simply by convergence of $\rho_k \to \rho$, since $G,K$ are bounded. Thus we have $$m \le J(\rho,u) \le \liminf_{k\to\infty} J(\rho_k,u_k) = m$$ proving that the pair $(\rho,u)$ is a solution of \eqref{eq:optProb}. \QED \\

Given this, we would like to prove first order necessary conditions for optimality akin to the classical Pontryagin Minimum Principle. To do so, we need a notions of an adjoint solution and a Hamiltonian. \\

For functions $\rho(x,t),q(x,t),u(x,t)$, we define the formal Lagrangian by \begin{equation} \label{eq:lagrangian}L(\rho,q,u) = J(\rho,u) + \int_{\Omega_T} q\bigg(\rho_t + \nabla \cdot (F(x,u)\rho)-\sum^{2}_{i,j=1} \frac{\partial^2}{\partial x_i \partial x_j}(D_{ij}(x,u)\rho) \bigg)dxdt \end{equation} Here the adjoint function $q(x,t)$ is playing the role of a Lagrange multiplier. Integrating by parts and expanding $J(\rho,u)$, we have 
\begin{align*}
L(\rho,q,u)&=\int_{\Omega_T} \ell(u(x,t))+ \rho\bigg(-q_t - F(x,u)\cdot \nabla q - \sum^2_{i,j=1} D_{i,j}(x,u) \frac{\partial^2q}{\partial x_i \partial x_j} + G(x,t)\bigg)dxdt \\
&\hspace{1cm} +\int_\Omega \rho(x,T)(q(x,T)+K(x)) - q(x,0)\rho_0(x)dx \\
&\hspace{1cm} + \int_0^T\int_{\partial \Omega} \rho(D\nabla q) \cdot n dS dt,
\end{align*} where we have enforced the zero flux boundary condition $E\cdot n = 0$. The adjoint equation is derived by enforcing  $$\bigg[\frac{d}{d\eta} L(\rho+\eta v, q, u)\bigg]_{\eta = 0} = 0$$ for arbitrary smooth test functions $v$ satisfying $v(x,0) = 0$. Computing formally, this yields \begin{equation}
    \label{eq:adjointEquation}\begin{aligned}
     &q_t + F(x,u)\cdot \nabla q + \sum^2_{i,j=1} D_{ij}(x,u)  \frac{\partial^2q}{\partial x_i \partial x_j} = G(x,t), \,\,\,\,\,\, \text{ in } \Omega \times [0,T), \\
     &q(x,T) = -K(x), \,\,\,\,\,\, \text{ in } \Omega,\\
     &(D\nabla q)\cdot n = 0, \,\,\,\,\,\,\,\,\,\,\,\,\, \text{ on } \partial \Omega \times [0,T).
    \end{aligned}
\end{equation} Solutions to \eqref{eq:adjointEquation} are considered in a weak sense analogous to \eqref{eq:weakSoln} above. That is, $q$ is a weak solution of \eqref{eq:adjointEquation} if $q(\cdot,T) = -K(\cdot)$ and \begin{equation}
    \label{eq:qpid}
    \int_\Omega \bigg(\psi q_t - q \nabla\cdot (F(x,u) \psi) - \sum^2_{i,j=1} \frac{\partial}{\partial x_i}(D_{ij}(x,u)\psi) \frac{\partial q}{\partial x_j} \bigg)dx= \int_\Omega G(x,t)\psi dx
\end{equation} for all $\psi \in H^1(\Omega)$ and a.e. $t \in [0,T]$. For fixed $u \in U_{\text{ad}}$, existence and uniqueness follow directly as in \cite[Ch. 7.1]{Evans}. Regularity can then be upgraded to $W^{2,p}(\Omega_T)$ via standard parabolic regularity results \cite[Thm. 7.35]{lieberman1996second}. This is the main reason that we need so much additional regularity for our control space. With controls uniformly bounded $W^{1,\infty}(\Omega_T)$ and taking values in the compact rectangle $R_U$, the maps $(x,t) \mapsto F(x,u(x,t))$ and $(x,t)\mapsto D(x,u(x,t))$ are in $W^{1,\infty}(\Omega_T)$ which allows for the increased regularity. \\

\thm If $K \in H^1(\Omega)$ and $G \in L^\infty(\Omega_T)$, then for any $u \in U_{\text{ad}},$ there is a unique weak solution $q$ of \eqref{eq:adjointEquation} with $$q \in L^2(0,T;H^1(\Omega)) \,\,\,\, \text{ and } \,\,\,\, q_t \in L^2(0,T; H^{-1}(\Omega)).$$ Further, if $K \in W^{2,p}(\Omega)$, then $q \in W^{2,1}_p(\Omega_T)$ and  $\|q\|_{W^{2,1}_p(\Omega_T)}$ is uniformly bounded over choices $u \in U_{\text{ad}}.$  \\

We henceforth assume that $K \in W^{2,\infty}(\Omega)$ so that we have uniformly bounded $\|q\|_{W^{2,1}_p(\Omega_T)}$ for any finite $p$. \\ 

We then formally define the Hamiltonian \begin{equation}
\label{eq:Hamiltonian}
H(\rho,q,u) = \ell(u) - \rho F(x,u)\cdot \nabla q + \sum^2_{i,j=1} \frac{\partial}{\partial x_j}(D_{ij}(x,u) \rho) \frac{\partial q}{\partial x_i}. 
\end{equation} Note that $H$ depends on $(x,t)$ implicitly through $\rho,q,u$; we have suppressed this dependence to reduce the notational burden. 

We prove that the optimal control map can be characterized as a pointwise minimizer of this Hamiltionian. Before doing so, we prove a lemma relating the cost functional and the Hamiltonian.\\

\lem For two control maps $u,\hat u \in U_{\text{ad}}$ with respective solutions $\rho, \hat \rho$ of \eqref{eq:controlledFP}, we have $$J(\rho,u) - J(\hat \rho, \hat u) = \int^T_0 \int_{\Omega} \Big[H(\hat \rho, q, u) - H(\hat \rho, q,\hat u) \Big] dx dt$$ where $q$ is the solution of \eqref{eq:adjointEquation} corresponding to $(\rho,u).$\\

\prf We compute \begin{equation} \label{eq:Jeq1}
\begin{aligned}
    J(\rho,u)-J(\hat \rho,\hat u) &= \int^T_{0}\int_\Omega \Big[\ell(u)-\ell(\hat u) + G(x,t)(\rho-\hat \rho)\Big] dxdt + \int_\Omega (\rho(\cdot,T)-\hat \rho(\cdot,T))K(x)dx\\
    &= \int^T_0\int_\Omega (\ell(u)-\ell(\hat u)) dxdt + \int_\Omega (\hat \rho(\cdot,T)-\rho(\cdot,T))q(\cdot, T)dx \\
    &\hspace{1cm} +\int^T_0 \int_\Omega \bigg((\rho - \hat \rho)q_t + (\rho - \hat \rho)F(x,u)\cdot \nabla q \bigg) dxdt \\&\hspace{1cm}-\int^T_0 \int_\Omega \sum^2_{i,j=1} \frac \partial{\partial x_j}(D_{ij}(x,u)[\rho - \hat \rho])  \frac{\partial q}{\partial x_i} dxdt\\
    &= \int^T_0\int_\Omega (\ell(u)-\ell(\hat u)) dxdt \\ 
    &\hspace{1cm} +\int^T_0 \int_\Omega \bigg((\hat \rho_t - \rho_t)q + (\rho - \hat \rho)F(x,u)\cdot \nabla q \bigg) dxdt \\&\hspace{1cm}-\int^T_0 \int_\Omega \sum^2_{i,j=1} \frac \partial{\partial x_j}(D_{ij}(x,u)[\rho - \hat \rho])  \frac{\partial q}{\partial x_i} dxdt\\
    \end{aligned}
\end{equation} where in the middle equality we have used \eqref{eq:qpid} with $\psi = \rho - \hat \rho,$ and in the last we have integrated by parts in $t$ and canceled boundary terms. But now with $q$ playing the role of the test function $\varphi$ in \eqref{eq:weakSoln}, we see that all terms with $\rho$ drop out of \eqref{eq:Jeq1} since $\rho$ is the weak solution of \eqref{eq:controlledFP} corresponding to $u$. Thus we arrive at \begin{equation}\label{eq:Jeq2}
\begin{aligned}
    J(\rho,u)-J(\hat \rho,\hat u) &= \int^T_0\int_\Omega \bigg(\ell(u)-\ell(\hat u) + \hat \rho_tq - \hat \rho F(x,u)\cdot \nabla q \bigg) dxdt\\ 
    &\hspace{1cm}+\int^T_0 \int_\Omega \sum^2_{i,j=1} \frac \partial{\partial x_j}(D_{ij}(x,u) \hat \rho)  \frac{\partial q}{\partial x_i}\bigg) dxdt\\
    &= \int^T_0 \int_\Omega \bigg(\ell(u) - \hat \rho F(x,u)\cdot q + \sum^2_{i,j=1} \frac \partial{\partial x_j}(D_{ij}(x,u) \hat \rho)  \frac{\partial q}{\partial x_i}\bigg) \\
    &\hspace{1cm}- \int^T_0 \int_\Omega \bigg(\ell(\hat u) - \hat \rho F(x,\hat u)\cdot q + \sum^2_{i,j=1} \frac \partial{\partial x_j}(D_{ij}(x,\hat u) \hat \rho)  \frac{\partial q}{\partial x_i}\bigg) \\ 
    &= \int^T_0 \int_\Omega \Big[H(\hat \rho, q,u) - H(\hat \rho,q,\hat u) \Big] dx dt
     \end{aligned} 
\end{equation} where in the last equality we have replaced $\int^T_0 \int_\Omega \hat \rho_t qdxdt$ using that $\hat \rho$ is a weak solution of \eqref{eq:controlledFP}. This completes the proof. \QED\\

\thm Suppose that the state-control pair $(\rho,u)$ is a solution of \eqref{eq:optProb} and that $q$ is the corresponding solution of \eqref{eq:adjointEquation}. Then for any $w \in R_U$, \begin{equation} \label{eq:Hineq}H(\rho(x,t),q(x,t),u(x,t)) \le H(\rho(x,t),q(x,t),w) \,\,\,\, \text{ for almost all } \,\, (x,t)\in \Omega_T.\end{equation} That is, $u(x,t)$ is a pointwise minimizer of the Hamiltonian when evaluated at the optimal state and costate pair. \\

\prf The proof utilizes the classical needle variation technique. Supposing that $(\rho,u)$ is optimal, we fix a point $(x_0,t_0) \in \Omega_T$ and $w \in R_u$. Define the set $$B_k(x_0,t_0) \defeq \{(x,t) \in \Omega_T \,\, : \,\, \abs{t-t_0} + \|x-x_0\|_2 \le 1/k \}$$ and the control map $$ u_k(x,t) = \begin{cases}
    u(x,t), & (x,t) \in \Omega_T \setminus B_k(x_0,t_0), \\
    w, & (x,t) \in \Omega_T \cap B_k(x_0,t_0).
\end{cases}.$$ In what follows, $w$ is the fixed constant value in $R_U$, while $u = u(x,t)$ depends on $(x,t)$. Assuming $\rho_k$ is the state map corresponding to $u_k$ and $q_k$ is the adjoint map corresponding to $(\rho_k, u_k)$,  \lref{3.7} tells us that \begin{equation}\label{eq:Jeq3}\begin{aligned}
   0 \le  \frac{J(\rho_k,u_k) - J(\rho, u)}{\abs{B_k(x_0,t_0)}} &=\frac 1{\abs{B_k(x_0,t_0)}} \int^T_0 \int_{\Omega} \Big[H( \rho, q_k, u_k) - H(\rho, q_k,u) \Big] dx dt\\
    &= \frac 1{\abs{B_k(x_0,t_0)}}\int_{B_k(x_0,t_0)}  \Big[H( \rho, q_k, w) - H( \rho, q_k,u) \Big] dx dt\\
    &= \frac 1{\abs{B_k(x_0,t_0)}}\int_{B_k(x_0,t_0)}  \Big[H( \rho, q, w) - H(\rho,q,u)\Big]dxdt \\&\hspace{1cm}  +\frac 1{\abs{B_k(x_0,t_0)}}\int_{B_k(x_0,t_0)}\Big[H(\rho,q_k,w) - H( \rho, q,w)\Big]dxdt\\ &\hspace{1cm}+\frac 1{\abs{B_k(x_0,t_0)}}\int_{B_k(x_0,t_0)} \Big[H(\rho,q,u) - H(\rho,q_k,u) \Big] dx dt.\\
    &=\frac 1{\abs{B_k(x_0,t_0)}}\int_{B_k(x_0,t_0)}  \Big[H( \rho, q, w) - H(\rho,q,u)\Big]dxdt + I_1 + I_2 
\end{aligned}\end{equation} where $I_1,I_2$ are the indented integrals. Now the map $$(x,t) \mapsto \Big[H(\rho(x,t),q(x,t),w) - H(\rho(x,t),q(x,t),u(x,t))\Big]$$ is in $L^1(\Omega_T)$, so by the Lebesgue differentiation theorem \begin{align*} \lim_{k\to\infty}\frac 1{\abs{B_k(x_0,t_0)}}\int_{B_k(x_0,t_0)}  &\Big[H( \rho, q, w) - H(\rho,q,u)\Big]dxdt \\ &= H(\rho(x_0,t_0),q(x_0,t_0),w)-H(\rho(x_0,t_0),q(x_0,t_0),u(x_0,t_0))\end{align*} as long as $(x_0,t_0) \in \Omega_T$ is a Lebesgue point for this function. Thus, taking the limit as $k\to\infty$ in \eqref{eq:Jeq3}, we arrive at \eqref{eq:Hineq} so long as $I_1,I_2 \to 0$. We see \begin{equation} \label{eq:I1} \begin{aligned}
    I_1 &= \frac 1 {\abs{B_k(x_0,t_0)}} \int_{B_k(x_0,t_0)} \Big[H(\rho,q_k,w) - H( \rho, q,w)\Big]dxdt \\
    &= \frac 1 {\abs{B_k(x_0,t_0)}} \int_{B_k(x_0,t_0)} \bigg[\rho F(x,w)\cdot [\nabla q - \nabla q_k] - \sum^2_{i,j=1} \frac{\partial}{\partial x_j}(D_{ij}(x,w)\rho)\frac{\partial}{\partial x_i}(q-q_k) \bigg]dxdt \\ 
    &= \frac 1 {\abs{B_k(x_0,t_0)}} \int_{B_k(x_0,t_0)} E(x,w)\cdot (\nabla q - \nabla q_k)dxdt\\
    &\le \left(\frac 1 {\abs{B_k(x_0,t_0)}} \int_{B_k(x_0,t_0)} \abs{E(x,w)}dxdt\right) \|\nabla q -\nabla q_k\|_{L^\infty(\Omega_T)}
\end{aligned} \end{equation} where $E$ is the Fokker-Planck flux defined in \eqref{eq:FLUX}. Now $q_k \to q$ strongly in $L^2(0,T; H^1(\Omega))$ by a standard
parabolic stability estimate \cite[Ch. 7, Thm. 2]{Evans}. However, we recall that for sufficiently regular $K$, $\|q\|_{W^{2,1}_p(\Omega_T)}$ is uniformly bounded over choices of controls for any $p > 1$. Taking $p$ large enough, we have that $W^{2,1}_p(\Omega_T)$ is compactly embedded in $C^{1,1/2}(\overline \Omega_T)$ \cite[Thm. 3.14]{Hu}, \cite[Ch. II, Lem. 3.3]{ladyz}. This implies that $\nabla q_k \to \nabla q$ uniformly along a subsequence. Assuming $(x_0,t_0)$ is a Lebesgue point of $\abs{E(x,w)}$, the term in the parentheses remains bounded as $k\to\infty$ and thus $I_1 \to 0$ by uniform convergence of $\nabla q_k$. The integral $I_2$ is treated identically. Because such points $(x_0,t_0)$ comprise almost all of $\Omega_T$, passing to the limit in \eqref{eq:Jeq3} shows that 
$$0\le H(\rho(x_0,t_0),q(x_0,t_0),w)-H(\rho(x_0,t_0),q(x_0,t_0),u(x_0,t_0)) \,\,\,\, \text{ for almost all } \,\, (x_0,t_0) \in \Omega_T$$ as desired. \QED

Given this, the first order necessary optimality conditions for this system are as follows. \\

\cor If $(\rho,u)$ is a solution of \eqref{eq:optProb}, then there exists a costate map $q$ such that \begin{equation}\label{eq:FP} \tag{FP}
\left\{\begin{aligned}
    &\rho_t + \nabla \cdot(F(x,u)\rho) = \sum_{i,j=1}^2 \frac{\partial^2}{\partial x_i \partial x_j} (D_{ij}(x,u)\rho), \,\,\,\,\, \text{ in } \Omega \times (0,T],\\
    &\rho(x,0) = \rho_0(x) \,\, \text{ in } \, \Omega, \,\,\,\, \,\,\,\,\, E(x,u) \cdot n= 0 \, \text{ on } \partial \Omega \text{ for } E \text{ as defined in } \eqref{eq:FLUX},
\end{aligned}\right.
\end{equation}
\begin{equation}
    \label{eq:ADJ}\tag{ADJ}
    \left\{\begin{aligned}
        &q_t + F(x,u)\cdot \nabla q + \sum^2_{i,j=1} D_{ij}(x,u)  \frac{\partial^2q}{\partial x_i \partial x_j} = G(x,t), \,\,\,\,\,\, \text{ in } \Omega \times [0,T), \hphantom{Hii} \\
     &q(x,T) = -K(x) \,\,\, \text{ in } \Omega, \,\,\,\,\,\, (D\nabla q)\cdot n = 0, \,\, \text{ on } \partial \Omega, 
    \end{aligned}\right.
\end{equation}
\vspace{0.25cm}
\begin{equation}
    \label{eq:OPT}\tag{OPT}
    u(x,t) \in \argmin_{w \in R_u} H(\rho(x,t),q(x,t),w) \,\,\,\,\, \text{ for a.e. } (x,t) \in \Omega \times [0,T].\hphantom{HHHH}
\end{equation}

\section{Numerical approximation}\label{sec:computation}

In this section, we describe the numerical methods used to solve the
FP optimal control problem \eqref{eq:optProb}. The numerical
procedure consists of finite-difference approximations of the forward
FP equation \eqref{eq:controlledFP} and the backward adjoint
equation \eqref{eq:adjointEquation}, coupled through the sequential
quadratic Hamiltonian (SQH) method. The latter iteratively updates the
space- and time-dependent feedback control
$u=(u_1,u_2,u_3)$ by pointwise minimization of a quadratically augmented
Hamiltonian.

\subsection{Numerical approximation of the forward FP equation}

Let
\[
x_i=a_x+(i-1)\Delta x,\qquad
y_j=a_y+(j-1)\Delta y,\qquad
t_n=n\Delta t,
\]
denote a uniform discretization of
$\Omega_T=\Omega\times(0,T)$, and let
$f_{i,j}^n$ denote the numerical approximation of the probability
density $\rho(x_i,y_j,t_n)$. Since $\rho$ represents a probability
density, the numerical treatment of the forward FP
equation should, as far as possible, retain the fundamental structural
properties of the continuous problem. In particular, the numerical
solution should remain nonnegative, conserve total probability under
the zero-flux boundary condition, and accurately balance drift and
diffusion without introducing spurious oscillations. For this reason,
we use a two-dimensional Chang--Cooper discretization for the drift
and diagonal diffusion terms.

We write the controlled FP equation as
\begin{equation}\label{eq:FPnumericalform}
\rho_t
=
-\partial_{x_1}(B_1\rho)
-\partial_{x_2}(B_2\rho)
+\partial_{x_1x_1}(D_{11}\rho)
+2\partial_{x_1x_2}(D_{12}\rho)
+\partial_{x_2x_2}(D_{22}\rho),
\end{equation}
where
\[
B_i=F_i(x,u),\qquad i=1,2.
\]
The diagonal contributions in \eqref{eq:FPnumericalform} can be written
in conservative flux form as
\[
-\partial_{x_1}
\left[
\left(B_1-\partial_{x_1}D_{11}\right)\rho
-
D_{11}\partial_{x_1}\rho
\right]
\]
and
\[
-\partial_{x_2}
\left[
\left(B_2-\partial_{x_2}D_{22}\right)\rho
-
D_{22}\partial_{x_2}\rho
\right].
\]
This motivates the modified drift coefficients
\[
\widetilde B_1
=
B_1-\partial_{x_1}D_{11},
\qquad
\widetilde B_2
=
B_2-\partial_{x_2}D_{22}.
\]
In the implementation, these quantities are evaluated at the
corresponding cell faces using finite differences of the diagonal
diffusion coefficients.

At an $x_1$-face $(x_{i+1/2},y_j)$, define
\[
w_{1,i+\frac12,j}^n
=
\frac{\Delta x\,\widetilde B_{1,i+\frac12,j}^n}
     {D_{11,i+\frac12,j}^n},
\]
and analogously at an $x_2$-face,
\[
w_{2,i,j+\frac12}^n
=
\frac{\Delta y\,\widetilde B_{2,i,j+\frac12}^n}
     {D_{22,i,j+\frac12}^n}.
\]
The Chang--Cooper weights are determined by
\begin{equation}\label{eq:CCdelta}
\delta(w)
=
\frac{1}{w}
-
\frac{1}{e^w-1},
\end{equation}
with the limiting value
\[
\delta(0)=\frac12.
\]
The quantity $w$ measures the relative strength of drift and diffusion
across a grid cell. Thus, unlike a standard centered discretization,
the Chang--Cooper weighting automatically biases the numerical flux in
a manner determined by the local drift--diffusion balance. For example, the discrete flux through the face
$(x_{i+1/2},y_j)$ may be written in the form
\[
\mathcal F_{i+\frac12,j}^{\,n}
=
\widetilde B_{1,i+\frac12,j}^{\,n}
\left[
(1-\delta_{1,i+\frac12,j}^{\,n})f_{i+1,j}^{\,n}
+
\delta_{1,i+\frac12,j}^{\,n}f_{i,j}^{\,n}
\right]
-
D_{11,i+\frac12,j}^{\,n}
\frac{f_{i+1,j}^{\,n}-f_{i,j}^{\,n}}{\Delta x},
\]
where
\[
\delta_{1,i+\frac12,j}^{\,n}
=
\delta\left(w_{1,i+\frac12,j}^n\right).
\]
Similarly,
\[
\mathcal G_{i,j+\frac12}^{\,n}
=
\widetilde B_{2,i,j+\frac12}^{\,n}
\left[
(1-\delta_{2,i,j+\frac12}^{\,n})f_{i,j+1}^{\,n}
+
\delta_{2,i,j+\frac12}^{\,n}f_{i,j}^{\,n}
\right]
-
D_{22,i,j+\frac12}^{\,n}
\frac{f_{i,j+1}^{\,n}-f_{i,j}^{\,n}}{\Delta y}.
\]
The diagonal drift--diffusion contribution at an interior grid point is
therefore represented by the conservative flux difference
\[
-\frac{
\mathcal F_{i+\frac12,j}^{\,n}
-
\mathcal F_{i-\frac12,j}^{\,n}}
{\Delta x}
-
\frac{
\mathcal G_{i,j+\frac12}^{\,n}
-
\mathcal G_{i,j-\frac12}^{\,n}}
{\Delta y}.
\]

A principal advantage of the Chang--Cooper construction is that the
drift and diffusion terms are not discretized independently. Instead,
the weighting in \eqref{eq:CCdelta} is chosen from their local ratio.
For one-dimensional FP equations, this construction is
well known to preserve nonnegativity under the appropriate time-step
restriction and to reproduce the correct discrete equilibrium
associated with zero probability flux. In the present two-dimensional
implementation, the same construction is applied direction-by-direction
to the diagonal drift--diffusion contributions. This significantly
reduces the artificial numerical diffusion and spurious oscillations
that may arise from a straightforward centered approximation.

The diffusion matrix in the present stochastic SIRS model is not
diagonal, so the FP equation also contains the mixed term
\[
2\partial_{x_1x_2}(D_{12}\rho).
\]
This contribution is treated separately. Setting
\[
g_{i,j}^n=D_{12,i,j}^n f_{i,j}^n,
\]
we use the centered approximation
\begin{equation}\label{eq:mixedFD}
\partial_{x_1x_2}(D_{12}\rho)_{i,j}^n
\approx
\frac{
g_{i+1,j+1}^n-g_{i+1,j-1}^n
-g_{i-1,j+1}^n+g_{i-1,j-1}^n
}
{4\Delta x\Delta y}.
\end{equation}
Consequently, if $\mathcal L_{\mathrm{CC}}^n$ denotes the
Chang--Cooper discretization of the diagonal drift--diffusion terms
and $\mathcal L_{\mathrm{mix}}^n$ denotes the approximation
\eqref{eq:mixedFD}, the full explicit update can be written schematically as
\begin{equation}\label{eq:FPexplicit}
f_{i,j}^{n+1}
=
f_{i,j}^n
+
\Delta t
\left[
(\mathcal L_{\mathrm{CC}}^n f^n)_{i,j}
+
2(\mathcal L_{\mathrm{mix}}^n f^n)_{i,j}
\right].
\end{equation}
Thus the diagonal part is advanced using Chang--Cooper fluxes, while
the cross-diffusion generated by the common transmission noise is
included through a centered mixed-derivative correction. This is the
same decomposition used in the numerical implementation. 

The time discretization in \eqref{eq:FPexplicit} is explicit, and
therefore the time step must be chosen sufficiently small relative to
the spatial mesh and the magnitudes of the drift and diffusion
coefficients. Under the corresponding stability restriction, the
Chang--Cooper component of the scheme retains the positivity-preserving
character of the underlying discretization. Because the mixed
derivative is added separately by a centered explicit approximation,
positivity of the complete two-dimensional update is not automatic for
an arbitrary time step. We therefore use a sufficiently fine temporal
grid and monitor the minimum value of the computed density throughout
the calculation.

An equally important property is conservation of total probability.
For the continuous FP equation,
\[
\frac{d}{dt}\int_\Omega \rho(x,t)\,dx
=
-\int_{\partial\Omega}E(x,u)\cdot n\,dS,
\]
and hence the zero-flux boundary condition implies
\[
\int_\Omega \rho(x,t)\,dx
=
\int_\Omega\rho_0(x)\,dx
=
1.
\]
The discrete Chang--Cooper terms are written as differences of
neighboring numerical fluxes, so the interior contributions telescope
when summed over all grid cells. Thus conservation reduces to the
treatment of the boundary fluxes.

Accordingly, after every forward Euler update, the zero-flux boundary
condition from \eqref{eq:controlledFP} is imposed directly. Because the
diffusion matrix contains the off-diagonal term $D_{12}$, the normal
flux cannot be enforced using only the diagonal diffusion coefficients.
For example, on a vertical boundary the condition involves
\[
E_1
=
B_1\rho
-
\partial_{x_1}(D_{11}\rho)
-
\partial_{x_2}(D_{12}\rho)
=
0,
\]
whereas on a horizontal boundary,
\[
E_2
=
B_2\rho
-
\partial_{x_1}(D_{12}\rho)
-
\partial_{x_2}(D_{22}\rho)
=
0.
\]
The boundary values of the density are therefore reconstructed so that
the complete normal probability flux, including the cross-diffusion
contribution, vanishes. The four corner values are obtained by
averaging the neighboring boundary values. With this treatment, the spatial discretization retains the
conservative character of the FP equation. 

\subsection{Numerical approximation of the adjoint equation}

The adjoint equation \eqref{eq:adjointEquation} is solved backward in
time from the terminal condition
\[
q(x,T)=-K(x).
\]
Unlike the forward equation, which is written in conservative
FP form, the adjoint equation is discretized directly in
its backward Kolmogorov form,
\[
q_t
+
F_1q_{x_1}
+
F_2q_{x_2}
+
D_{11}q_{x_1x_1}
+
2D_{12}q_{x_1x_2}
+
D_{22}q_{x_2x_2}
=
G.
\]
The face-centered drift coefficients used by the forward
Chang--Cooper solver are first averaged to the spatial grid points.
For interior grid points, the first derivatives are approximated by
centered differences,
\[
(q_{x_1})_{i,j}
\approx
\frac{q_{i+1,j}-q_{i-1,j}}{2\Delta x},
\qquad
(q_{x_2})_{i,j}
\approx
\frac{q_{i,j+1}-q_{i,j-1}}{2\Delta y},
\]
while the diagonal second derivatives are approximated by the usual
three-point formulas. The mixed derivative is approximated by
\[
(q_{x_1x_2})_{i,j}
\approx
\frac{
q_{i+1,j+1}-q_{i+1,j-1}
-q_{i-1,j+1}+q_{i-1,j-1}}
{4\Delta x\Delta y}.
\]
Consequently, if $\mathcal L_h$ denotes the resulting spatial
finite-difference operator, the backward time update takes the form
\[
q^{n}
=
q^{n+1}
+
\Delta t\left(\mathcal L_h q^{n+1}-G^{n+1}\right).
\]
The adjoint boundary condition
\[
(D\nabla q)\cdot n=0
\]
is also enforced directly at every backward time step. Thus, on the
vertical boundaries the discretization accounts for the coupling
between $D_{11}q_{x_1}$ and $D_{12}q_{x_2}$, while on the horizontal
boundaries it accounts for the coupling between
$D_{12}q_{x_1}$ and $D_{22}q_{x_2}$. The corner values are again
obtained by averaging adjacent boundary values.

\subsection{A sequential quadratic Hamiltonian method}

We employ the sequential quadratic Hamiltonian (SQH) method to
numerically resolve the optimal feedback maps. A full exposition of the SQH method with worked examples can be found in \cite{SQH}. 

The SQH method is an iterative procedure based on successive quadratic perturbations of the
Hamiltonian \eqref{eq:Hamiltonian}. For a penalization parameter
$\varepsilon>0$ and a current control approximation $\widetilde u$, we
define the augmented Hamiltonian
\begin{equation}\label{eq:augmented_Hamiltonian}
H_\varepsilon(\rho,q,u,\widetilde u)
=
H(\rho,q,u)
+
\varepsilon |u-\widetilde u|_2^2 .
\end{equation}
The quadratic term restricts the magnitude of the change in the
control between consecutive SQH iterations. In particular, increasing
$\varepsilon$ forces the new control to remain closer to the previous
iterate, whereas decreasing $\varepsilon$ permits a larger update.

Starting from an initial control $u^0$, the corresponding forward
FP solution $\rho^0$ and adjoint solution $q^0$ are
computed. At iteration $k$, the new feedback control is obtained
pointwise in $(x,t)$ by solving
\begin{equation}\label{eq:SQHmin}
u^{k+1}(x,t)
\in
\arg\min_{w\in R_U}
H_\varepsilon
\left(
\rho^k(x,t),q^k(x,t),w,u^k(x,t)
\right).
\end{equation}
This pointwise update is particularly convenient for the present
problem because the three control components can be obtained
explicitly and subsequently projected onto their admissible intervals.

For the mixed $L^1/L^2$ control cost used here, the pointwise
minimization leads to explicit candidate values for each control.
Writing
\[
A=q_{x_1x_1}-2q_{x_1x_2}+q_{x_2x_2},
\]
the updates implemented numerically are
\[
u_1^{k+1}
=
\Pi_{[u_{1,\min},u_{1,\max}]}
\left(
\frac{
-\beta_1+2\varepsilon u_1^k
+\beta\rho SI(q_{x_1}-q_{x_2})
-2\sigma_\beta\rho S^2I^2A
}{
\beta_2+2\varepsilon
-2\sigma_\beta\rho S^2I^2A
}
\right),
\]
\[
u_2^{k+1}
=
\Pi_{[u_{2,\min},u_{2,\max}]}
\left(
\frac{
-\beta_1+2\varepsilon u_2^k
-\lambda\rho(1-S-I)q_{x_1}
-2\sigma_\lambda\rho(1-S-I)^2q_{x_1x_1}
}{
\beta_2+2\varepsilon
-2\sigma_\lambda\rho(1-S-I)^2q_{x_1x_1}
}
\right),
\]
and
\[
u_3^{k+1}
=
\Pi_{[u_{3,\min},u_{3,\max}]}
\left(
\frac{
-\beta_1+2\varepsilon u_3^k
-\gamma\rho Iq_{x_2}
+2\sigma_\gamma\rho I^2q_{x_2x_2}
}{
\beta_2+2\varepsilon
-2\sigma_\gamma\rho I^2q_{x_2x_2}
}
\right),
\]
where $\Pi_{[a,b]}$ denotes projection onto the interval $[a,b]$.
After updating the controls, the drift and diffusion coefficients are
recomputed and the forward FP equation is solved again.

The magnitude of the control update is measured by
\begin{equation}\label{eq:SQHupdate}
\tau_k
=
\|u^{k+1}-u^k\|_{(L^2(\Omega_T))^3}^2
=
\sum_{i=1}^3
\int_0^T\int_\Omega
|u_i^{k+1}-u_i^k|^2\,dx\,dt.
\end{equation}
The iteration terminates when $\tau_k<\kappa$ or when a prescribed
maximum number of iterations has been reached.

To determine whether the candidate control is accepted, we employ a
sufficient-decrease criterion. Specifically, the candidate is accepted
when
\begin{equation}\label{eq:SQHaccept}
J(\rho^{k+1},u^{k+1})-J(\rho^k,u^k)
\leq
-\eta\,\tau_k .
\end{equation}
If this condition fails, the candidate control and state are rejected,
the previous iterate is restored, and the quadratic penalization
parameter is increased according to
\[
\varepsilon\leftarrow \sigma\varepsilon,
\qquad \sigma>1.
\]
If the sufficient-decrease condition is satisfied, the candidate
control is accepted, the corresponding adjoint equation is solved, and
the penalization parameter is decreased according to
\[
\varepsilon\leftarrow\zeta\varepsilon,
\qquad 0<\zeta<1.
\]
Thus, the SQH procedure adaptively balances the size of the control
update against the decrease obtained in the objective functional.

\subsection{Analysis of the SQH method}

We next establish some properties of the SQH iteration. The first result
gives a sufficient-decrease property analogous to that obtained for
time-dependent controls, with the control norm now taken over the full
space-time domain $\Omega_T$.

\begin{theorem}\label{theoremSQH}
Let $(\rho^k,u^k)$ and $(\rho^{k+1},u^{k+1})$ be generated by the SQH
method applied to \eqref{eq:optProb}, with
$u^k,u^{k+1}\in U_{\mathrm{ad}}$. Then, for the current value
$\varepsilon>0$ of the SQH penalization parameter,
\begin{equation}\label{eq:decreaseJ}
J(\rho^{k+1},u^{k+1})
-
J(\rho^k,u^k)
\leq
-\varepsilon
\|u^{k+1}-u^k\|_{(L^2(\Omega_T))^3}^2.
\end{equation}
In particular, if $\varepsilon\geq\eta$, then
\[
J(\rho^{k+1},u^{k+1})
-
J(\rho^k,u^k)
\leq
-\eta\,\tau_k,
\]
where
\[
\tau_k
=
\|u^{k+1}-u^k\|_{(L^2(\Omega_T))^3}^2.
\]
\end{theorem}

\begin{proof}
By definition of the SQH update, for almost every
$(x,t)\in\Omega_T$,
\[
H_\varepsilon
(\rho^k,q^k,u^{k+1},u^k)
\leq
H_\varepsilon
(\rho^k,q^k,w,u^k)
\qquad
\text{for all }w\in R_U.
\]
In particular, choosing $w=u^k$ gives
\[
H_\varepsilon
(\rho^k,q^k,u^{k+1},u^k)
\leq
H_\varepsilon
(\rho^k,q^k,u^k,u^k)
=
H(\rho^k,q^k,u^k).
\]
Using the definition of the augmented Hamiltonian, we therefore obtain
\[
H(\rho^k,q^k,u^{k+1})
+
\varepsilon |u^{k+1}-u^k|_2^2
\leq
H(\rho^k,q^k,u^k)
\]
for almost every $(x,t)\in\Omega_T$. Integrating over $\Omega_T$ and
applying the Hamiltonian identity derived in Section~\ref{sec:analysis}
gives
\[
\begin{aligned}
J(\rho^{k+1},u^{k+1})
-
J(\rho^k,u^k)
&=
\int_0^T\int_\Omega
\Big[
H(\rho^k,q^k,u^{k+1})
-
H(\rho^k,q^k,u^k)
\Big]\,dx\,dt\\
&\leq
-\varepsilon
\|u^{k+1}-u^k\|_{(L^2(\Omega_T))^3}^2,
\end{aligned}
\]
which proves \eqref{eq:decreaseJ}.
\end{proof}

As a consequence of Theorem~\ref{theoremSQH}, we obtain a basic
convergence property of the accepted SQH iterates.

\begin{theorem}\label{theoremSQHterminates}
Suppose that the accepted SQH iterations satisfy
$\varepsilon_k\geq\eta>0$. Then
\begin{itemize}
\item[(a)] the sequence
$\{J(\rho^k,u^k)\}$ is decreasing and converges to some
\[
J^*\geq \inf_{\rho,u}J(\rho,u);
\]
\item[(b)]
\[
\lim_{k\to\infty}
\|u^{k+1}-u^k\|_{(L^2(\Omega_T))^3}
=0.
\]
\end{itemize}
\end{theorem}

\begin{proof}
Statement (a) follows immediately from
\eqref{eq:decreaseJ}, since the objective functional is bounded below.
For statement (b), using $\varepsilon_k\geq\eta$ gives
\[
\|u^{k+1}-u^k\|_{(L^2(\Omega_T))^3}^2
\leq
\frac{1}{\eta}
\left[
J(\rho^k,u^k)
-
J(\rho^{k+1},u^{k+1})
\right].
\]
Summing from $k=0$ to $K$ yields the telescoping estimate
\[
0
\leq
\sum_{k=0}^{K}
\|u^{k+1}-u^k\|_{(L^2(\Omega_T))^3}^2
\leq
\frac{1}{\eta}
\left[
J(\rho^0,u^0)
-
J(\rho^{K+1},u^{K+1})
\right].
\]
Since the objective functional is bounded below, the right-hand side
remains bounded as $K\to\infty$. Hence
\[
\sum_{k=0}^{\infty}
\|u^{k+1}-u^k\|_{(L^2(\Omega_T))^3}^2
<\infty,
\]
which implies
\[
\|u^{k+1}-u^k\|_{(L^2(\Omega_T))^3}
\longrightarrow0.
\]
\end{proof}

In particular, for any prescribed tolerance $\kappa>0$, the stopping
criterion
\[
\tau_k
=
\|u^{k+1}-u^k\|_{(L^2(\Omega_T))^3}^2
<\kappa
\]
is eventually reached along the accepted SQH iterates. Moreover, the
continuity properties of the control-to-state map established in
Section~\ref{sec:analysis} imply that the corresponding changes in the
FP states also tend to zero.

We finally give a stronger convergence statement for the mixed
$L^2$--$L^1$ control cost used in the numerical experiments. Define the
SQH update map
\[
\mathcal T_\varepsilon:
(L^2(\Omega_T))^3
\longrightarrow
(L^2(\Omega_T))^3
\]
by
\[
\mathcal T_\varepsilon(u)
=
\argmin_{w\in R_U}
H_\varepsilon
\big(
\rho_u,q_u,w,u
\big),
\]
where $\rho_u=\mathcal S(u)$ and $q_u$ denotes the corresponding
adjoint state.

\begin{theorem}[Convergence of SQH for the mixed $L^2$--$L^1$ cost]
\label{thm:SQHconvergence}
Assume that
\[
\ell(u)
=
\frac{\beta_2}{2}|u|_2^2
+
\beta_1|u|_1,
\qquad
\beta_1,\beta_2>0,
\]
and suppose that, for a fixed $\varepsilon>0$, the SQH update map
$\mathcal T_\varepsilon$ is Lipschitz continuous and satisfies
\begin{equation}\label{eq:SQHcontraction}
\|
\mathcal T_\varepsilon(u)
-
\mathcal T_\varepsilon(\widetilde u)
\|_{(L^2(\Omega_T))^3}
\leq
\theta
\|
u-\widetilde u
\|_{(L^2(\Omega_T))^3}
\end{equation}
for all admissible controls $u,\widetilde u$, with some
$0\leq\theta<1$. Then the SQH iteration
\[
u^{k+1}=\mathcal T_\varepsilon(u^k)
\]
converges strongly in $(L^2(\Omega_T))^3$ to a unique fixed point
$u^*$. Moreover, if
\[
\rho^*=\mathcal S(u^*)
\]
and $q^*$ is the corresponding adjoint state, then
\[
u^*(x,t)
\in
\argmin_{w\in R_U}
H(\rho^*(x,t),q^*(x,t),w)
\]
for almost every $(x,t)\in\Omega_T$.
\end{theorem}

\begin{proof}
Under assumption \eqref{eq:SQHcontraction},
$\mathcal T_\varepsilon$ is a contraction on
$(L^2(\Omega_T))^3$. Hence, by the Banach fixed-point theorem, there
exists a unique $u^*$ such that
\[
\mathcal T_\varepsilon(u^*)=u^*,
\]
and
\[
u^k\to u^*
\qquad\text{strongly in }(L^2(\Omega_T))^3.
\]
At the fixed point,
\[
u^*
=
\argmin_{w\in R_U}
H_\varepsilon(\rho^*,q^*,w,u^*).
\]
Since the quadratic SQH penalty vanishes at $w=u^*$,
\[
H_\varepsilon(\rho^*,q^*,u^*,u^*)
=
H(\rho^*,q^*,u^*).
\]
Therefore the fixed-point condition gives
\[
H(\rho^*,q^*,u^*)
\leq
H(\rho^*,q^*,w)
\qquad
\text{for all }w\in R_U
\]
for almost every $(x,t)\in\Omega_T$, which is precisely the
pointwise Hamiltonian minimum condition.
\end{proof}

\section{Numerical results}
\label{sec:results}

In this section, we present numerical results for the FP
optimal control framework computed using the SQH
algorithm. We set $b=\delta=0.01$, so that the total population in the
underlying SIRS model approaches unity. The state variables are the
susceptible and infected fractions, $x=(S,I)$, and the FP
equation is solved on the computational domain
$\Omega=[0,1]\times[0,1]$. We discretize each spatial direction using
$50$ grid points and take the final time to be $T=2$, with $10001$
temporal grid points.

We take the transmission, recovery, and loss-of-immunity parameters to be
$\beta=3$, $\gamma=1$, and $\lambda=0.5$, respectively. The three
feedback controls $u_1(S,I,t)$, $u_2(S,I,t)$, and $u_3(S,I,t)$ act by
reducing transmission, reducing loss of immunity, and increasing
recovery, respectively. Their admissible ranges are
\[
0\leq u_1\leq 1,\qquad
0\leq u_2\leq 1,\qquad
0\leq u_3\leq 0.5.
\]
Thus, unlike time-dependent open-loop interventions, the controls
obtained from the FP formulation are feedback policies whose
values depend on both the current epidemic state $(S,I)$ and time.

The initial probability density is taken to be a two-dimensional Gaussian
centered at $(S_0,I_0)=(0.99,0.01)$. Specifically, the numerical
implementation uses
\[
f_0(S,I)
=
\frac{1}{\pi d}
\exp\left(
-\frac{(S-0.99)^2+(I-0.01)^2}{d}
\right),
\qquad d=0.1,
\]
restricted to the computational domain and normalized numerically to have
unit mass. This represents uncertainty in the initial epidemic state
around a population consisting predominantly of susceptible individuals
with a small initially infected fraction.

Stochasticity is introduced through uncertainty in the transmission,
loss-of-immunity, and recovery mechanisms. In the controlled stochastic
system, the corresponding diffusion amplitudes are modulated by the same
feedback controls that modify the drift, so that the transmission noise is
scaled by $(1-u_1)$, the loss-of-immunity noise by $(1-u_2)$, and the
recovery noise by $(1+u_3)$. In the simulations presented here, the
baseline noise amplitudes are taken to be
$\sigma_\beta=\sigma_\lambda=\sigma_\gamma=0.5$.

The objective functional combines the state-dependent epidemic cost with
a mixed $L^1/L^2$ penalty on the three feedback controls. In particular,
the control contribution has the form
\[
\beta_1\int_0^T\int_\Omega
\left(|u_1|+|u_2|+|u_3|\right)\,dx\,dt
+
\frac{\beta_2}{2}\int_0^T\int_\Omega
\left(u_1^2+u_2^2+u_3^2\right)\,dx\,dt,
\]
where we take $\beta_1=0.02, \beta_2=0.05.$ The $L^1$ contribution penalizes the overall amount of intervention and
can promote sparse or concentrated control actions, while the quadratic
term regularizes the controls and discourages excessively large
interventions. \\

\begin{figure}[h!]
\centering

\includegraphics[width=0.31\textwidth]{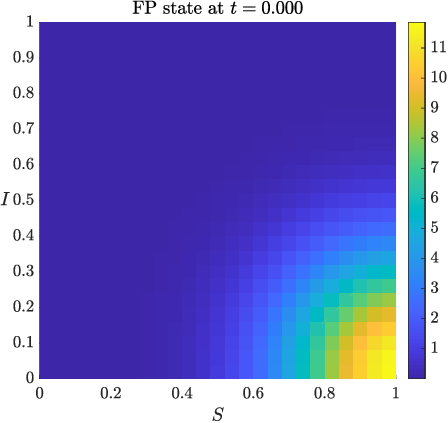}
\includegraphics[width=0.31\textwidth]{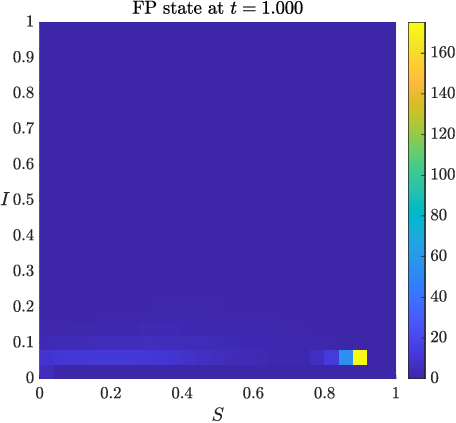}
\includegraphics[width=0.31\textwidth]{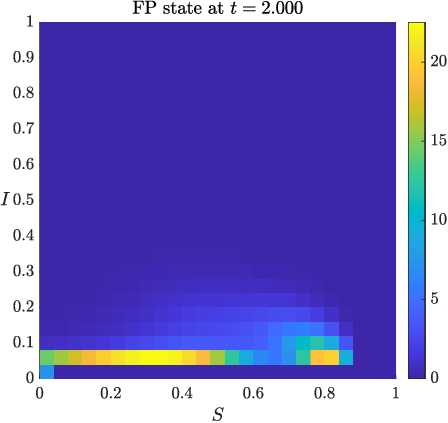}\\

\includegraphics[width=0.31\textwidth]{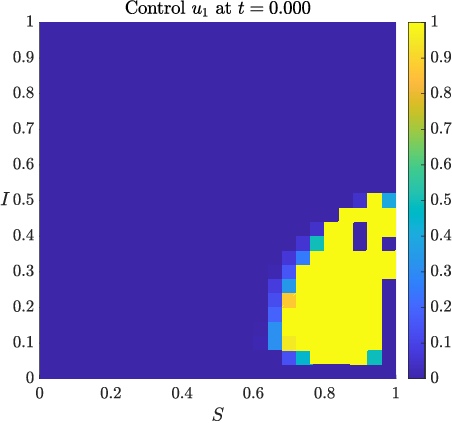}
\includegraphics[width=0.31\textwidth]{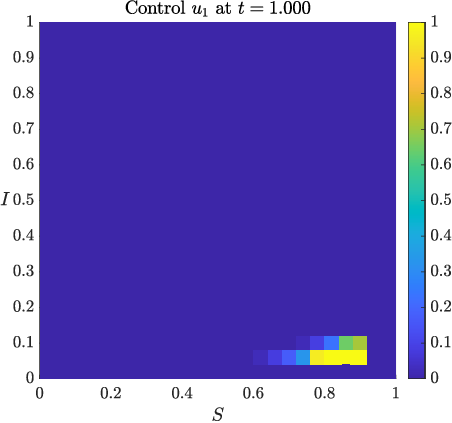}
\includegraphics[width=0.31\textwidth]{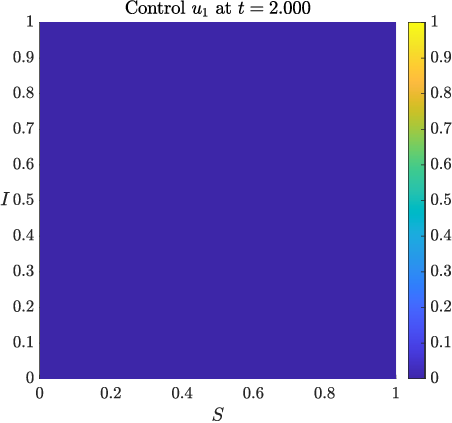}\\

\includegraphics[width=0.31\textwidth]{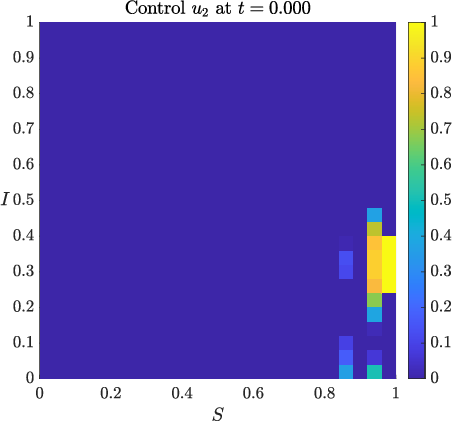}
\includegraphics[width=0.31\textwidth]{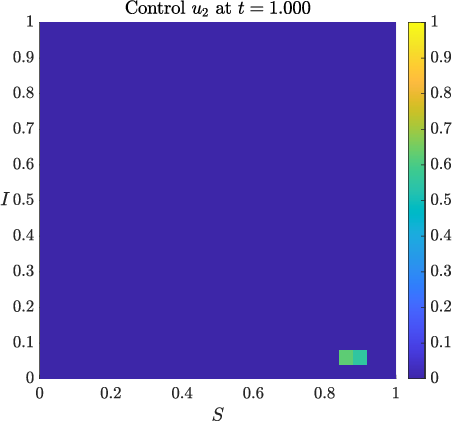}
\includegraphics[width=0.31\textwidth]{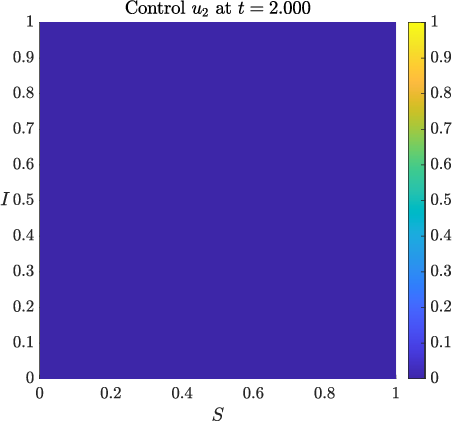}\\

\includegraphics[width=0.31\textwidth]{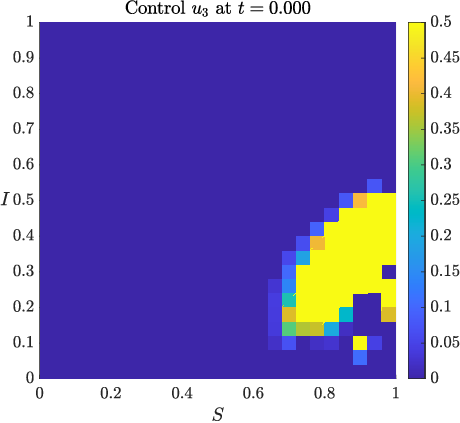}
\includegraphics[width=0.31\textwidth]{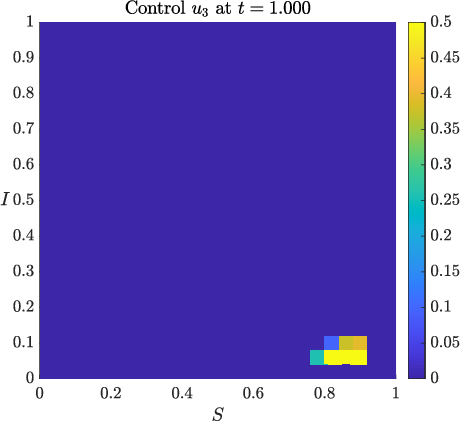}
\includegraphics[width=0.31\textwidth]{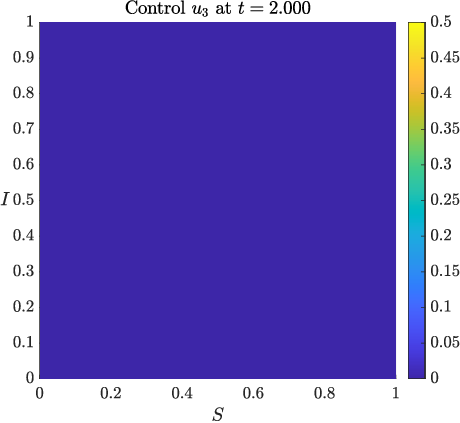}

\caption{Controlled FP solution and optimal feedback controls.
The columns correspond, from left to right, to the initial time $t=0$,
the midpoint $t=T/2$, and the final time $t=T$. The first row shows the
controlled probability density $f(S,I,t)$, while the second, third, and
fourth rows show the optimal feedback controls $u_1(S,I,t)$,
$u_2(S,I,t)$, and $u_3(S,I,t)$, respectively.}
\label{fig:controlled_FP_controls_Case1}
\end{figure}

\noindent {\bf Scenario 1.}
In the first scenario, we consider the running cost $G(S,I,t)=I,~ K(S,I)=0.$
Thus, infection is penalized continuously throughout the time interval,
with larger values of the infected fraction incurring a proportionally
larger state cost. No terminal-state cost is imposed. This provides a
baseline test of the feedback-control framework in which the objective
is to reduce the cumulative infection burden while balancing this
reduction against the cost of the three interventions. The resulting controlled FP density and optimal feedback
controls are shown in
Figure~\ref{fig:controlled_FP_controls_Case1}. The controlled density evolves toward states with smaller values of $I$,
while the feedback maps show how the three interventions are deployed
over the $(S,I)$ state space and time. The control $u_1$ reduces
transmission and hence suppresses the transition from susceptible to
infected individuals, $u_2$ reduces the loss of immunity and therefore
limits the return of recovered individuals to the susceptible class,
and $u_3$ increases the recovery rate and promotes the transition from
infected to recovered individuals. The distinct spatial and temporal
structures of these controls illustrate how the optimizer combines
prevention of new infections, retention of immunity, and enhanced
recovery to reduce the cumulative infection burden.

\begin{figure}[h!]
\centering

\includegraphics[width=0.31\textwidth]{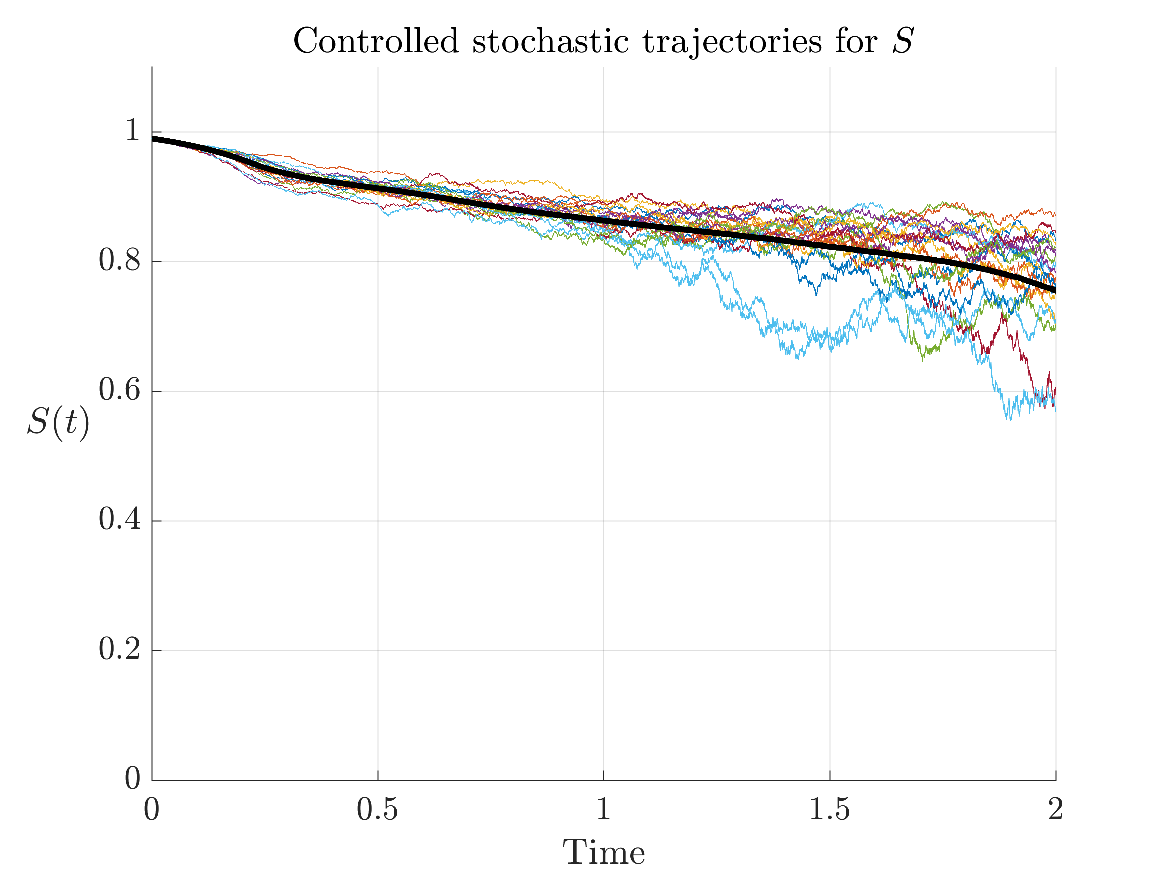}
\includegraphics[width=0.31\textwidth]{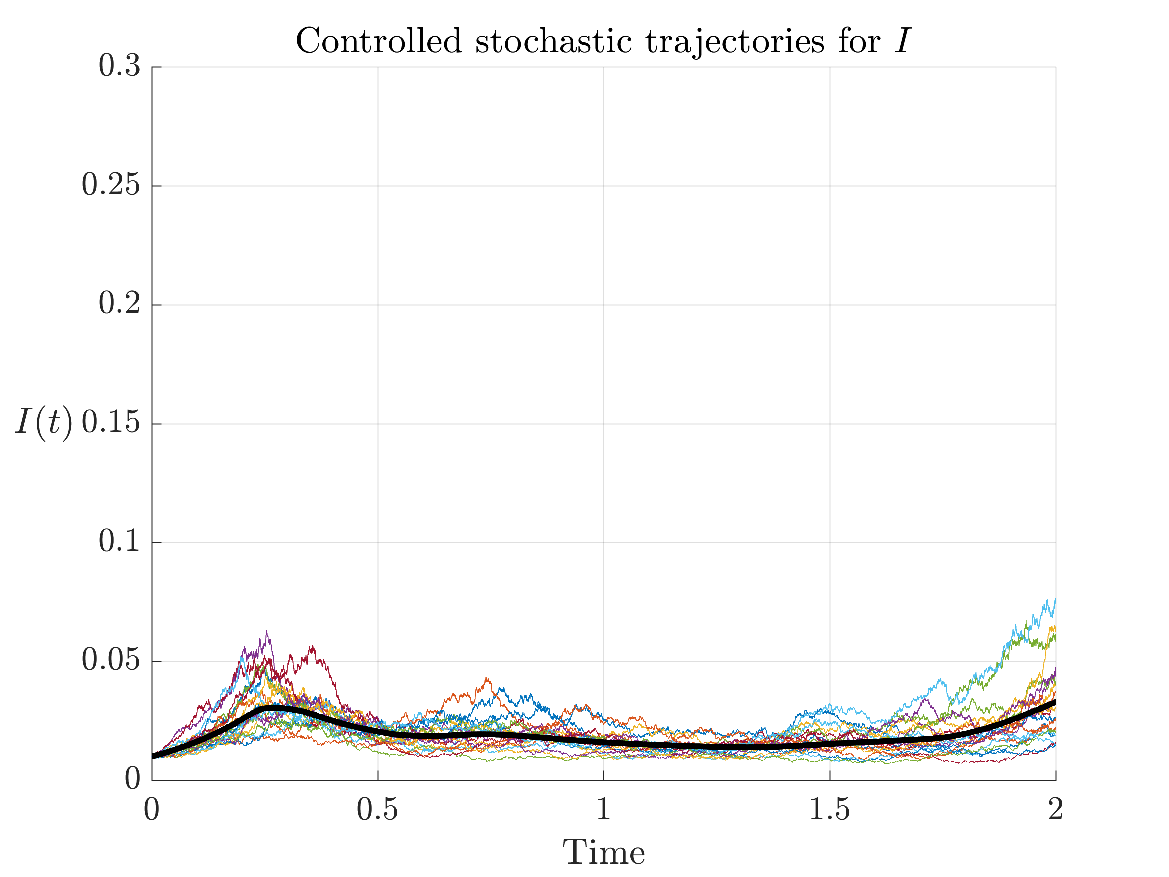}
\includegraphics[width=0.27\textwidth]{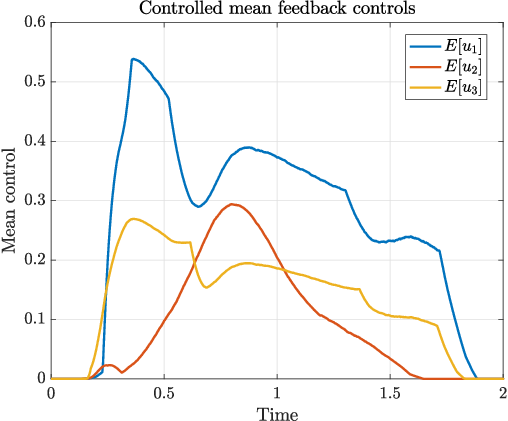}\\[-1mm]

\includegraphics[width=0.31\textwidth]{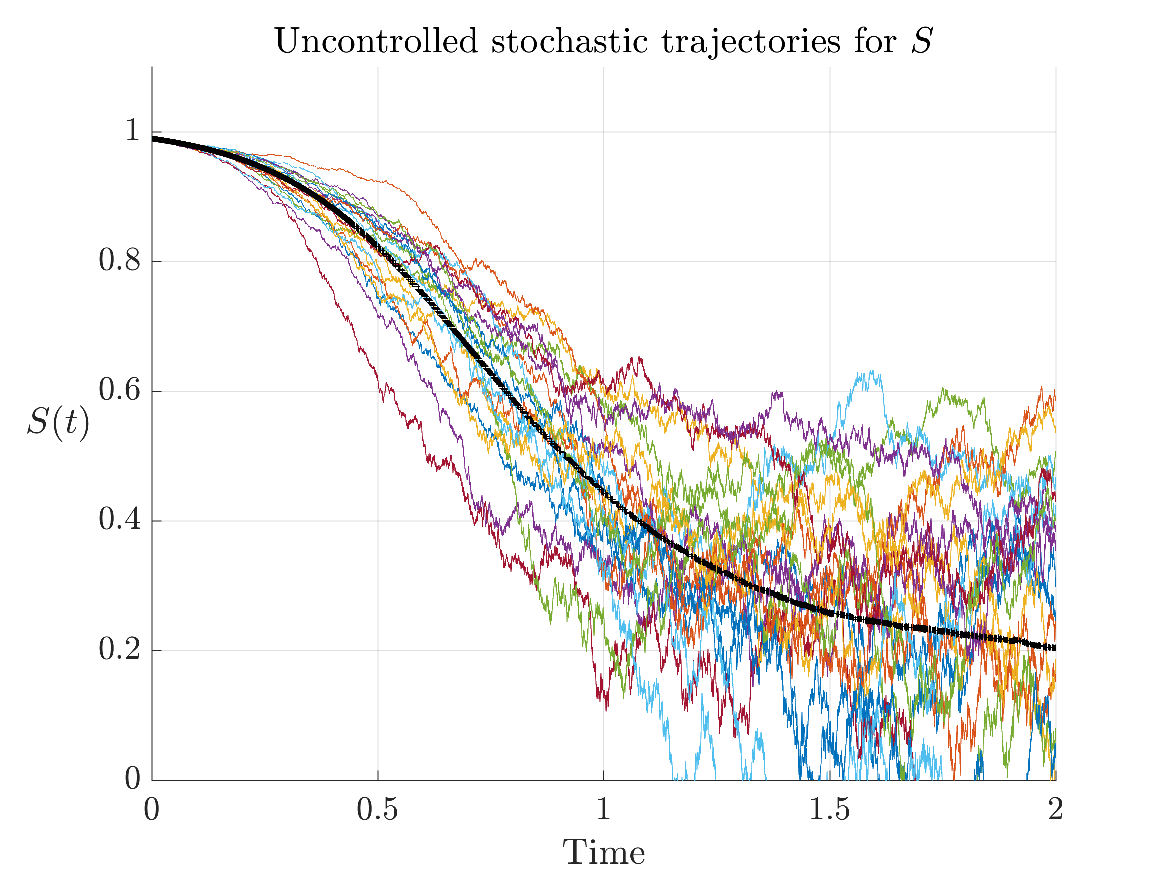}
\includegraphics[width=0.31\textwidth]{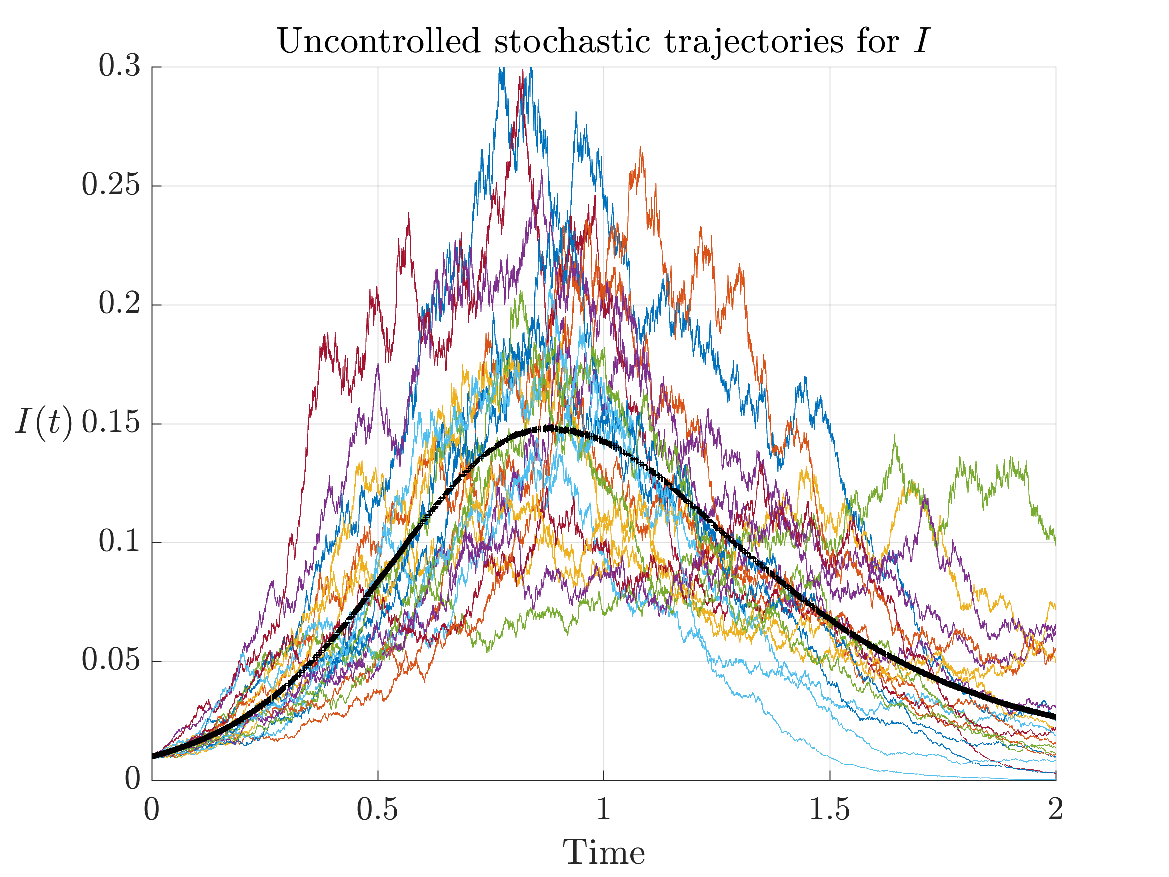}
\includegraphics[width=0.27\textwidth]{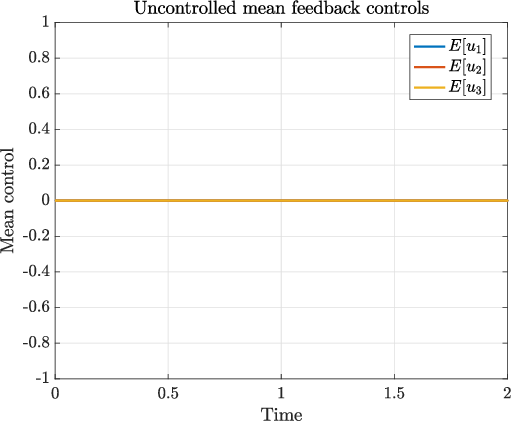}

\caption{Monte Carlo simulations of the controlled and uncontrolled
stochastic system. The top row shows the controlled stochastic
trajectories for the susceptible state $S$ and infected state $I$,
together with the mean feedback controls evaluated along the stochastic
trajectories. The bottom row shows the corresponding uncontrolled
stochastic trajectories and mean controls. The individual curves
represent stochastic realizations, while the emphasized curve denotes
the ensemble mean.}
\label{fig:controlled_uncontrolled_SDE_Case1}
\end{figure}

The corresponding controlled and uncontrolled stochastic simulations
are shown in
Figure~\ref{fig:controlled_uncontrolled_SDE_Case1}.  In the uncontrolled
system, the susceptible fraction decreases substantially while the
infected fraction undergoes a pronounced increase, whereas the optimized
feedback policy considerably limits both the loss of susceptible
individuals and the growth of the infected population. The ensemble
means exhibit the same qualitative behavior as the individual
realizations. The mean controls evaluated along the controlled
stochastic trajectories show how the state-dependent feedback policy is
implemented dynamically: $u_1$ provides the principal mechanism for
suppressing new infections, while $u_2$ acts by retaining individuals
in the recovered class and $u_3$ promotes removal from the infected
class through enhanced recovery. Thus, the stochastic simulations
confirm that the feedback controls reduce the infection burden relative
to the uncontrolled dynamics while adapting the intervention to the
evolving epidemic state.\\

\begin{figure}[h!]
\centering

\includegraphics[width=0.31\textwidth]{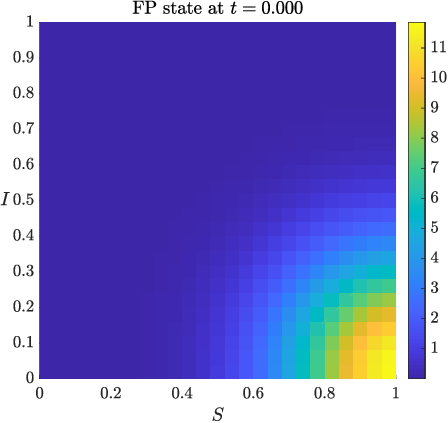}
\includegraphics[width=0.31\textwidth]{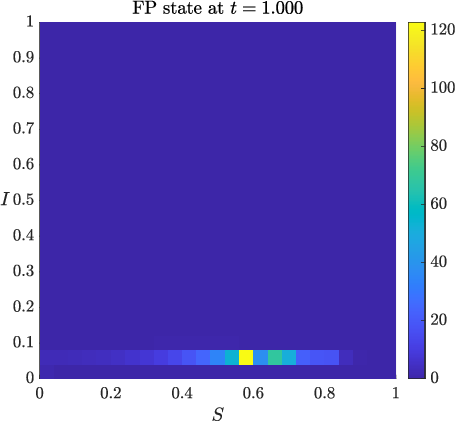}
\includegraphics[width=0.31\textwidth]{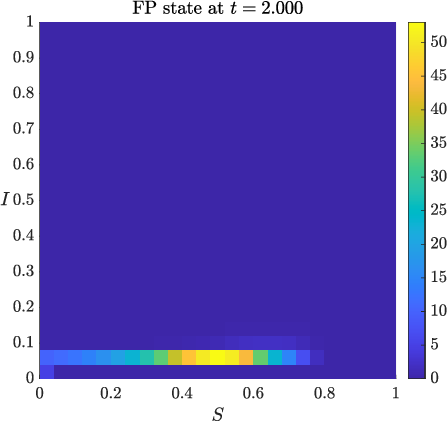}\\

\includegraphics[width=0.31\textwidth]{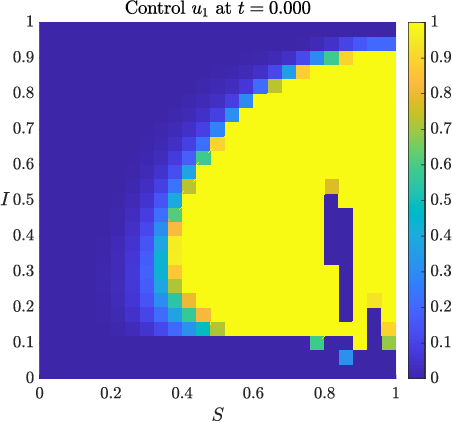}
\includegraphics[width=0.31\textwidth]{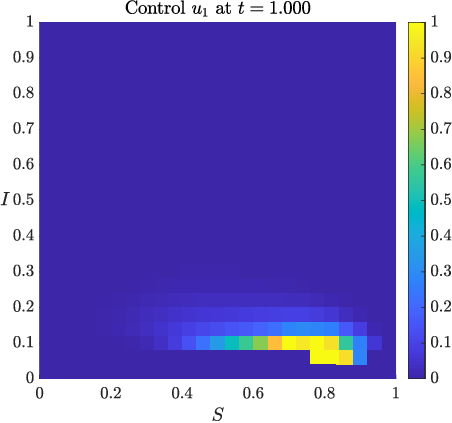}
\includegraphics[width=0.31\textwidth]{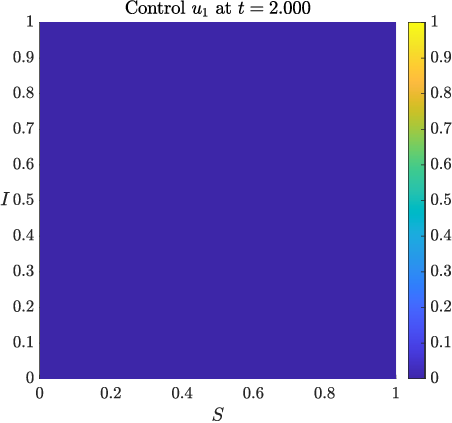}\\

\includegraphics[width=0.31\textwidth]{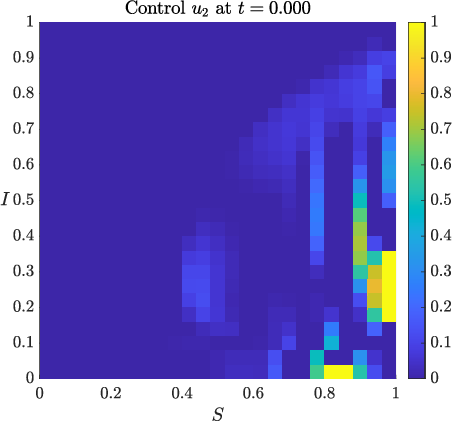}
\includegraphics[width=0.31\textwidth]{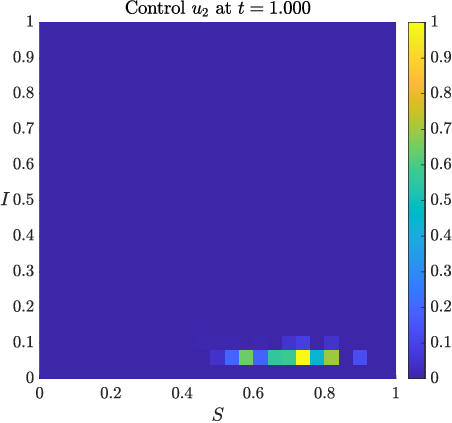}
\includegraphics[width=0.31\textwidth]{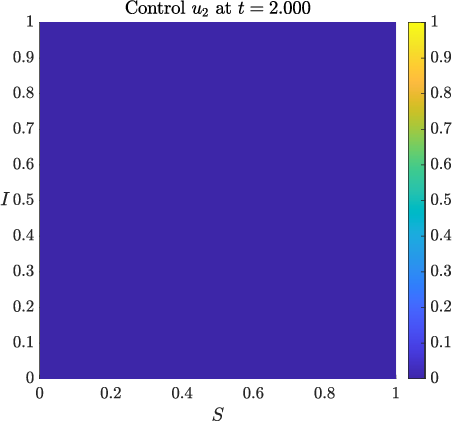}\\

\includegraphics[width=0.31\textwidth]{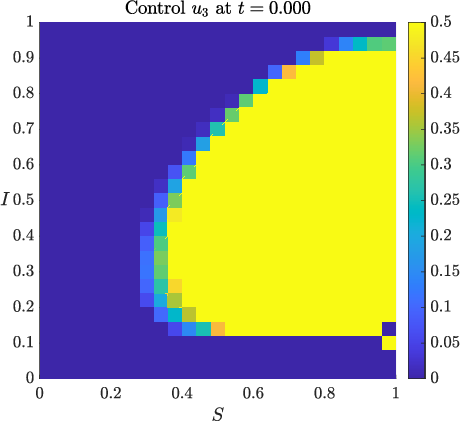}
\includegraphics[width=0.31\textwidth]{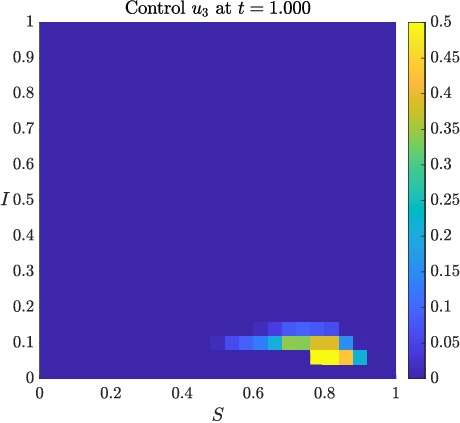}
\includegraphics[width=0.31\textwidth]{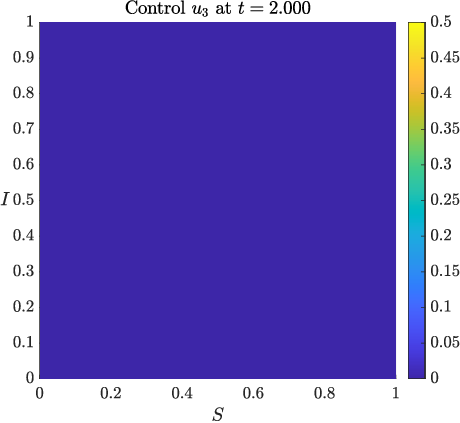}

\caption{Controlled FP solution and optimal feedback controls.
The columns correspond, from left to right, to the initial time $t=0$,
the midpoint $t=T/2$, and the final time $t=T$. The first row shows the
controlled probability density $f(S,I,t)$, while the second, third, and
fourth rows show the optimal feedback controls $u_1(S,I,t)$,
$u_2(S,I,t)$, and $u_3(S,I,t)$, respectively.}
\label{fig:controlled_FP_controls_Case2}
\end{figure}

\noindent {\bf Scenario 2.}
In this scenario, we consider running cost function
$
G(S,I,t)=100\,\mathbf{1}_{\{I\geq 0.08\}}(I-0.08),
$
 and zero terminal cost. This imposes no state penalty while the infected fraction remains below
$0.08$, but strongly penalizes excursions above this threshold. This is meant to represent the goal of ``flattening the curve" so as to not overwhelm the healthcare system. The controlled FP density and optimal feedback controls are
shown in Figure~\ref{fig:controlled_FP_controls_Case2}.  In contrast to
Scenario~1, the optimizer is not rewarded for reducing $I$ throughout
the entire state space; instead, control is valuable primarily when it
reduces the probability of entering or remaining in the penalized region
$I>0.08$. This is reflected in the state-dependent structure of the
feedback maps. The control $u_1$ reduces transmission and therefore
limits the flow from the susceptible to the infected class, $u_2$
reduces loss of immunity and retains individuals in the recovered class,
and $u_3$ increases recovery and promotes removal from the infected
class. Together, these mechanisms reshape the probability density so
as to limit its movement into states with infection levels above the
prescribed threshold, while avoiding unnecessary intervention in
regions where the infection-dependent running cost vanishes. 

\begin{figure}[h!]
\centering

\includegraphics[width=0.31\textwidth]{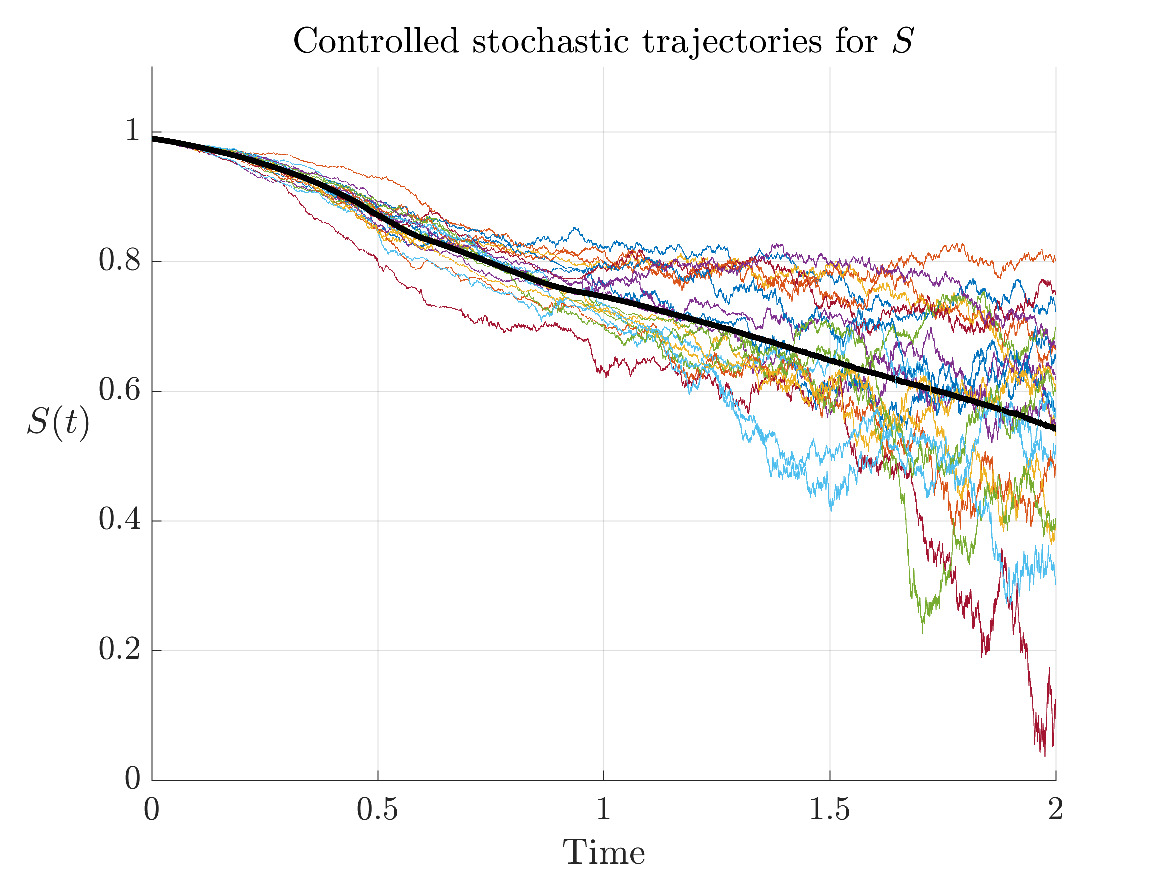}
\includegraphics[width=0.31\textwidth]{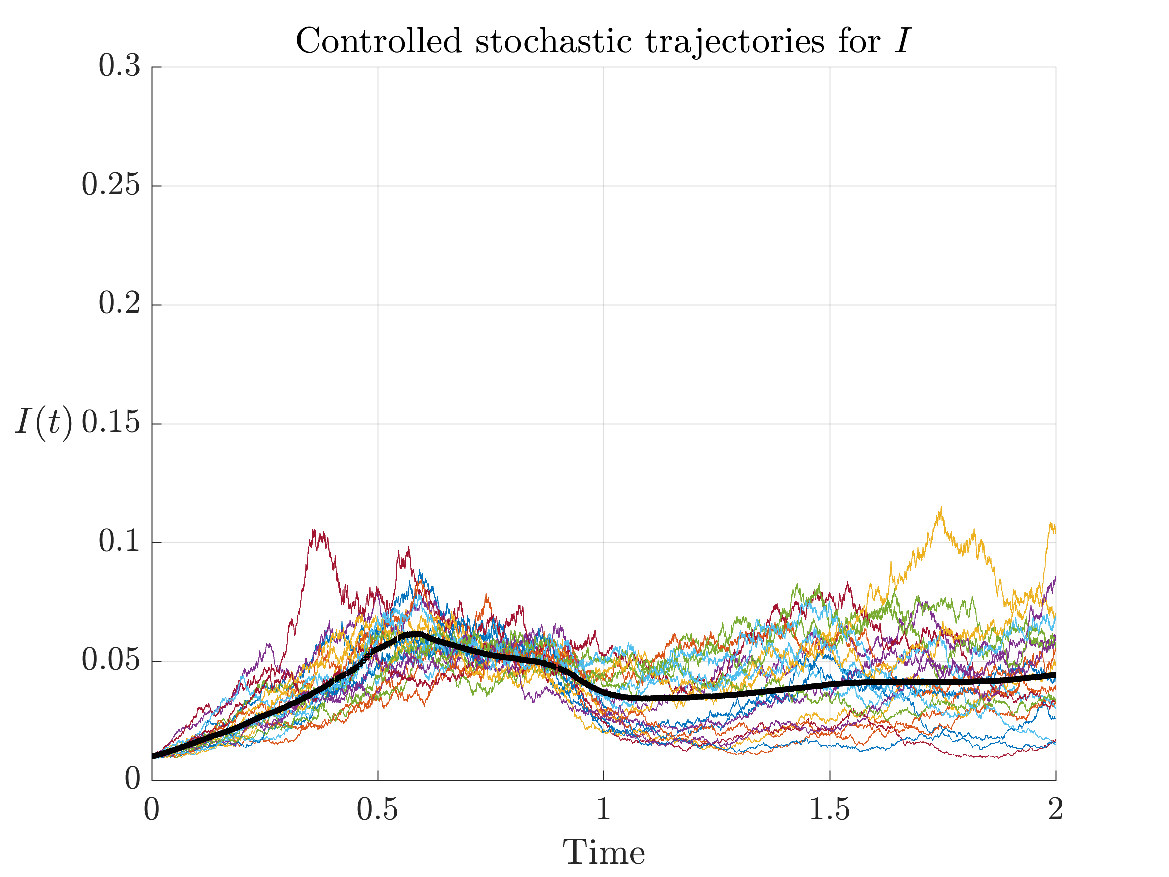}
\includegraphics[width=0.27\textwidth]{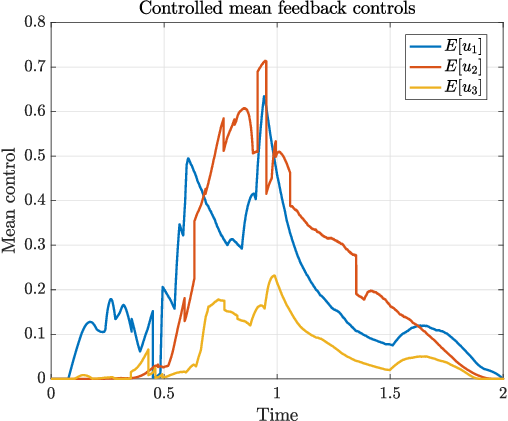}\\[-1mm]

\includegraphics[width=0.31\textwidth]{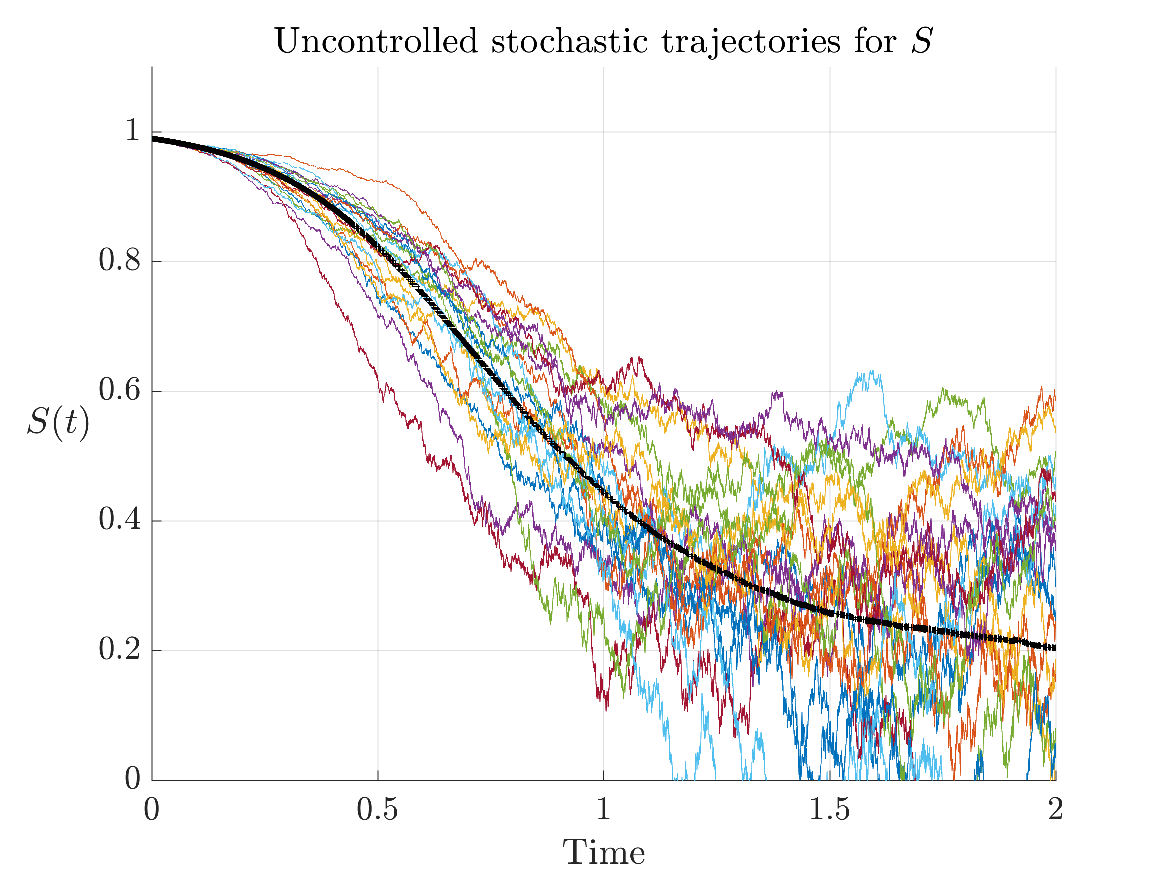}
\includegraphics[width=0.31\textwidth]{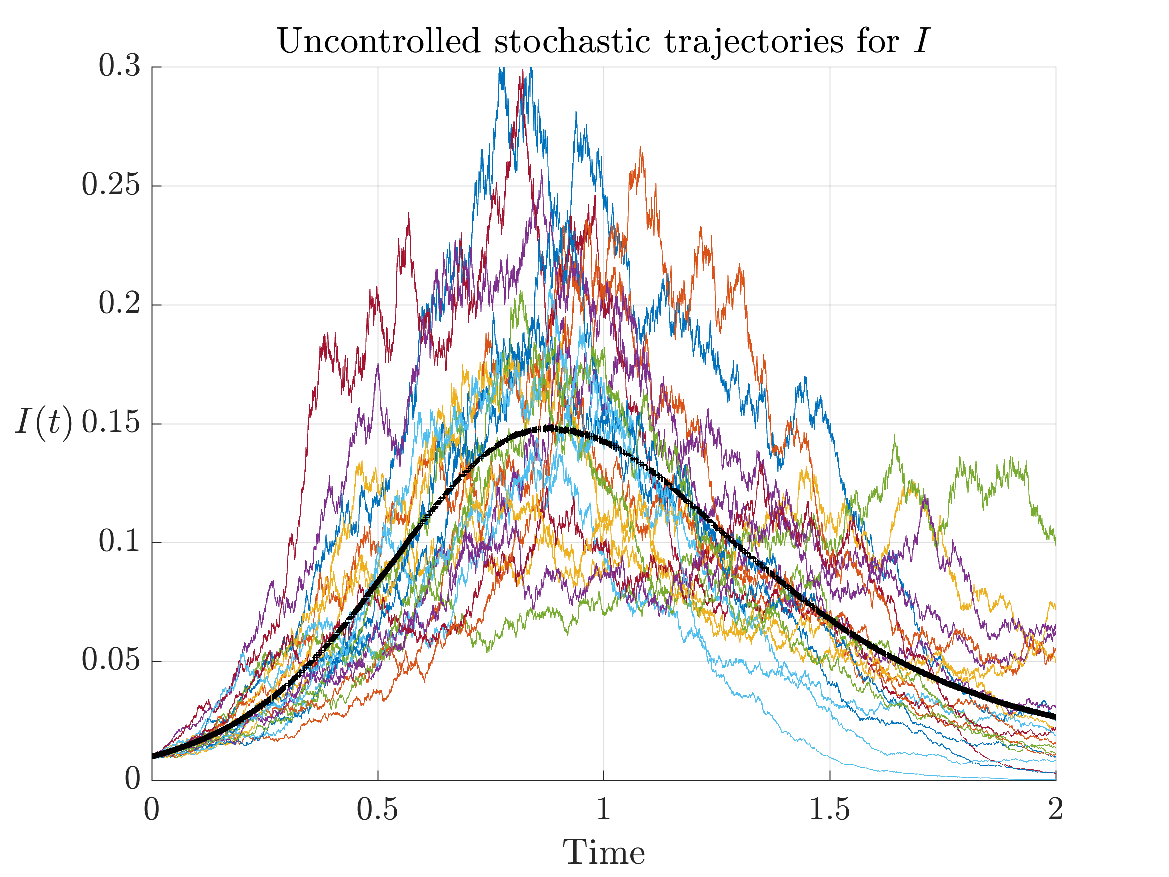}
\includegraphics[width=0.27\textwidth]{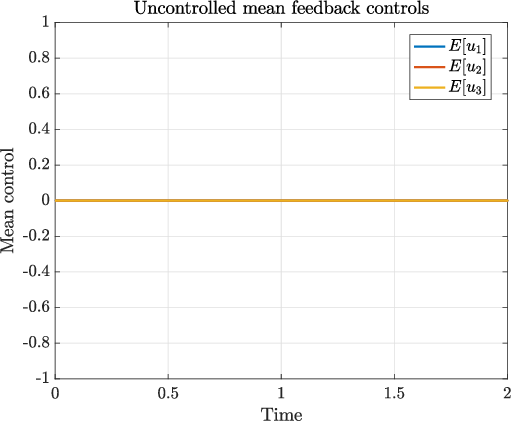}

\caption{Monte Carlo simulations of the controlled and uncontrolled
stochastic system. The top row shows the controlled stochastic
trajectories for the susceptible state $S$ and infected state $I$,
together with the mean feedback controls evaluated along the stochastic
trajectories. The bottom row shows the corresponding uncontrolled
stochastic trajectories and mean controls. The individual curves
represent stochastic realizations, while the emphasized curve denotes
the ensemble mean.}
\label{fig:controlled_uncontrolled_SDE_Case2}
\end{figure}

The corresponding controlled and uncontrolled stochastic simulations
are shown in
Figure~\ref{fig:controlled_uncontrolled_SDE_Case2}. 
The uncontrolled
trajectories exhibit a pronounced increase in the infected fraction,
with trajectories entering the region $I>0.08$ that is strongly
penalized by the objective, together with a substantial decrease in the
susceptible population. Under the optimized feedback policy, the growth
of the infected population is markedly suppressed and the stochastic
trajectories are largely kept away from the high-infection region,
while the decline in the susceptible population is also reduced. The
mean controls along the stochastic trajectories illustrate the
mechanism by which this threshold regulation is achieved: $u_1$
suppresses the generation of new infections, $u_2$ reduces the return
of recovered individuals to the susceptible class, and $u_3$ enhances
the removal of infected individuals through recovery. The resulting
behavior highlights the principal feature of this scenario: rather
than minimizing infection uniformly, the feedback strategy allocates
control effort specifically to reduce the probability and extent of
excursions above the prescribed infection threshold.\\

\begin{figure}[h!]
\centering

\includegraphics[width=0.31\textwidth]{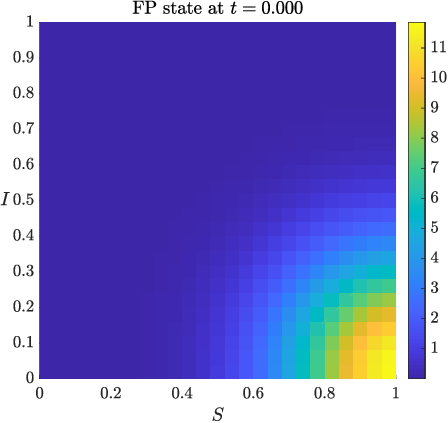}
\includegraphics[width=0.31\textwidth]{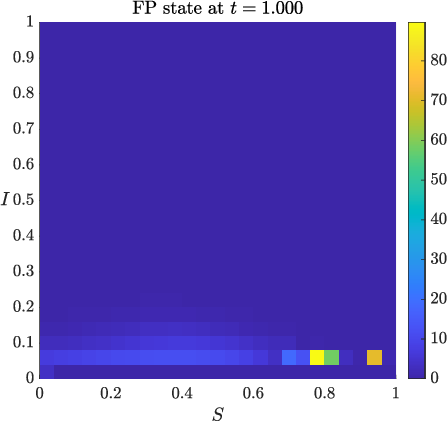}
\includegraphics[width=0.31\textwidth]{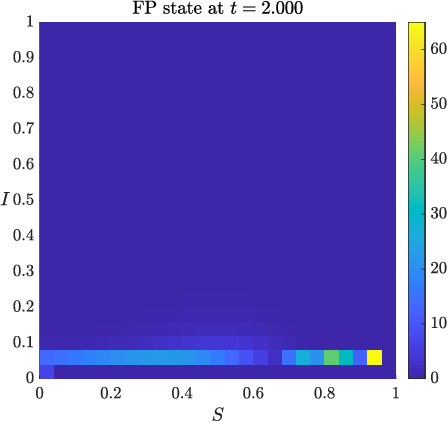}\\

\includegraphics[width=0.31\textwidth]{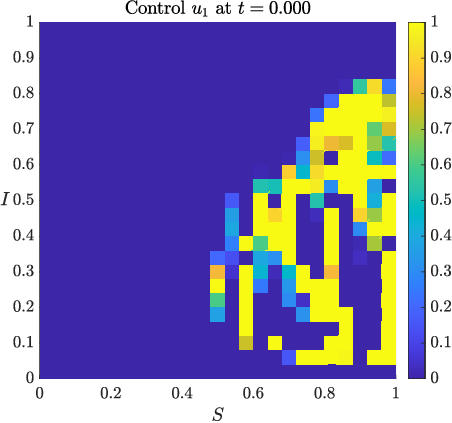}
\includegraphics[width=0.31\textwidth]{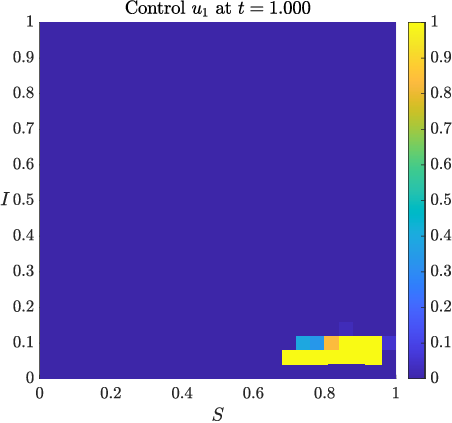}
\includegraphics[width=0.31\textwidth]{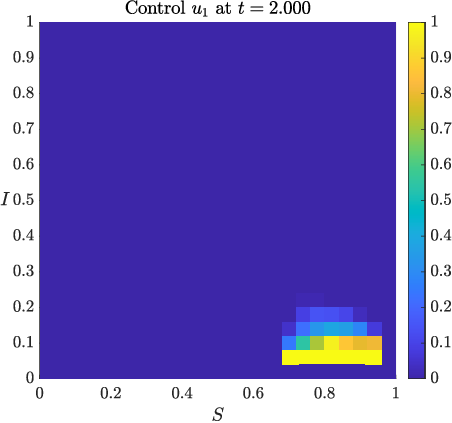}\\

\includegraphics[width=0.31\textwidth]{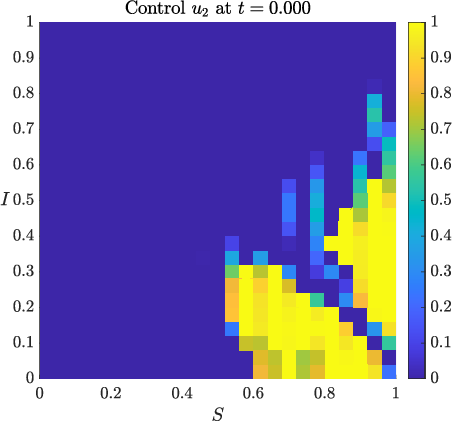}
\includegraphics[width=0.31\textwidth]{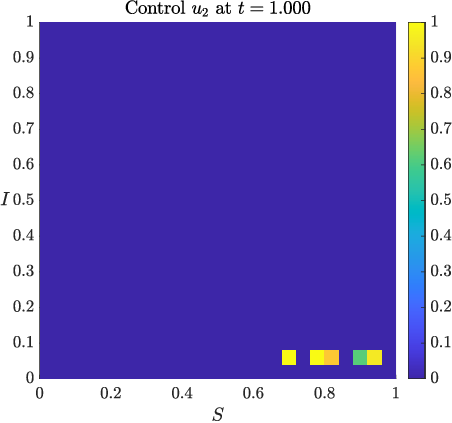}
\includegraphics[width=0.31\textwidth]{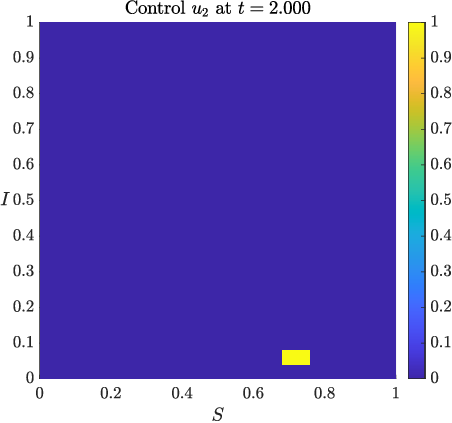}\\

\includegraphics[width=0.31\textwidth]{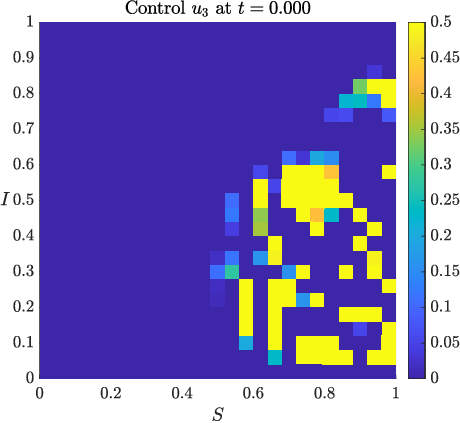}
\includegraphics[width=0.31\textwidth]{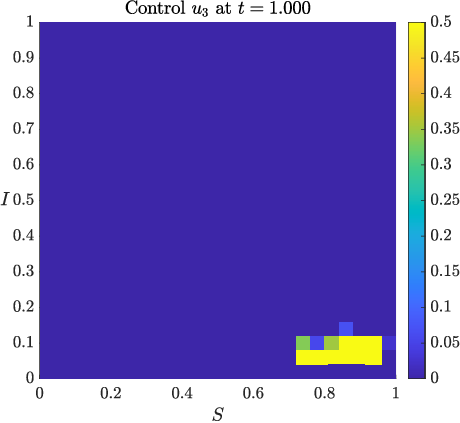}
\includegraphics[width=0.31\textwidth]{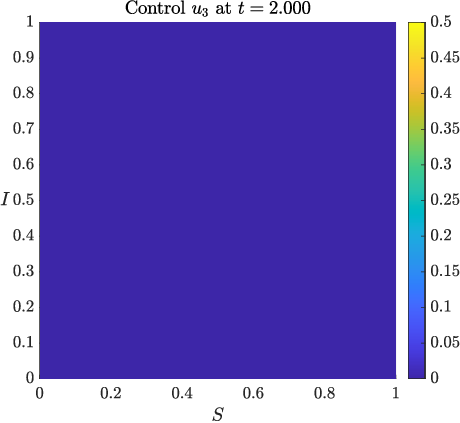}

\caption{Controlled FP solution and optimal feedback controls.
The columns correspond, from left to right, to the initial time $t=0$,
the midpoint $t=T/2$, and the final time $t=T$. The first row shows the
controlled probability density $f(S,I,t)$, while the second, third, and
fourth rows show the optimal feedback controls $u_1(S,I,t)$,
$u_2(S,I,t)$, and $u_3(S,I,t)$, respectively.}
\label{fig:controlled_FP_controls_Case3}
\end{figure}

\noindent {\bf Scenario 3.}
In this scenario, we consider
$
G(S,I,t)=0,~
K(S,I)=-10\max(S-\overline S,0),
$
so that there is no running state cost and the objective instead
provides a terminal reward for states in which the susceptible fraction
exceeds the prescribed level $\overline S$.The
resulting controlled FP solution and optimal feedback controls
are shown in Figure~\ref{fig:controlled_FP_controls_Case3}.  Since the
objective acts only through the terminal susceptible population, the
feedback controls are selected according to their ability to steer the
probability density toward states satisfying $S>\overline S$ at the final time.
The control $u_1$ contributes directly to this objective by reducing
transmission and hence limiting the loss of susceptible individuals to
the infected class. The control $u_2$ has a different effect: it reduces
the return of recovered individuals to the susceptible class and
therefore does not directly promote a larger value of $S$, while $u_3$
increases recovery from the infected to the recovered class and can
influence the subsequent evolution of both $I$ and $S$. The resulting
feedback maps therefore reflect the different contributions of these
three mechanisms to the terminal susceptible-state objective rather
than a direct penalty on the infected population during the time
interval.

\begin{figure}[h!]
\centering

\includegraphics[width=0.31\textwidth]{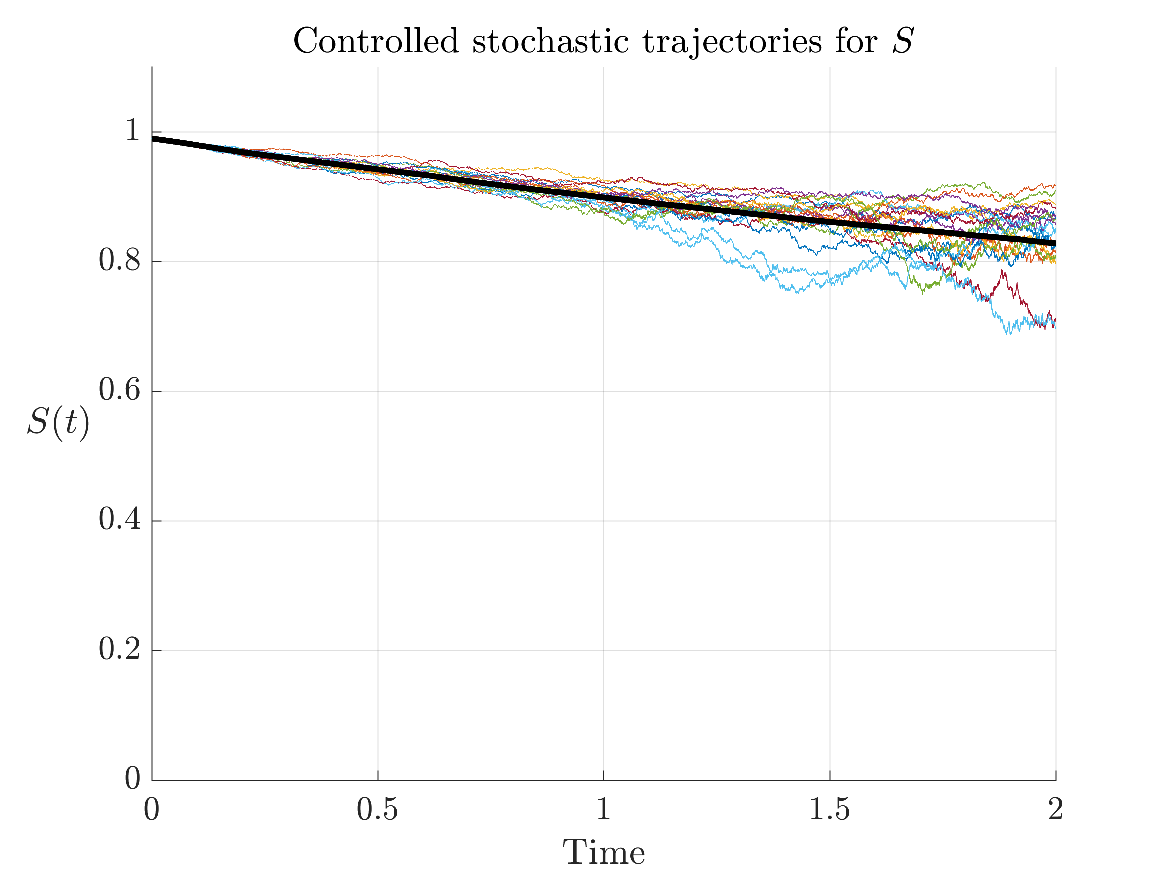}
\includegraphics[width=0.31\textwidth]{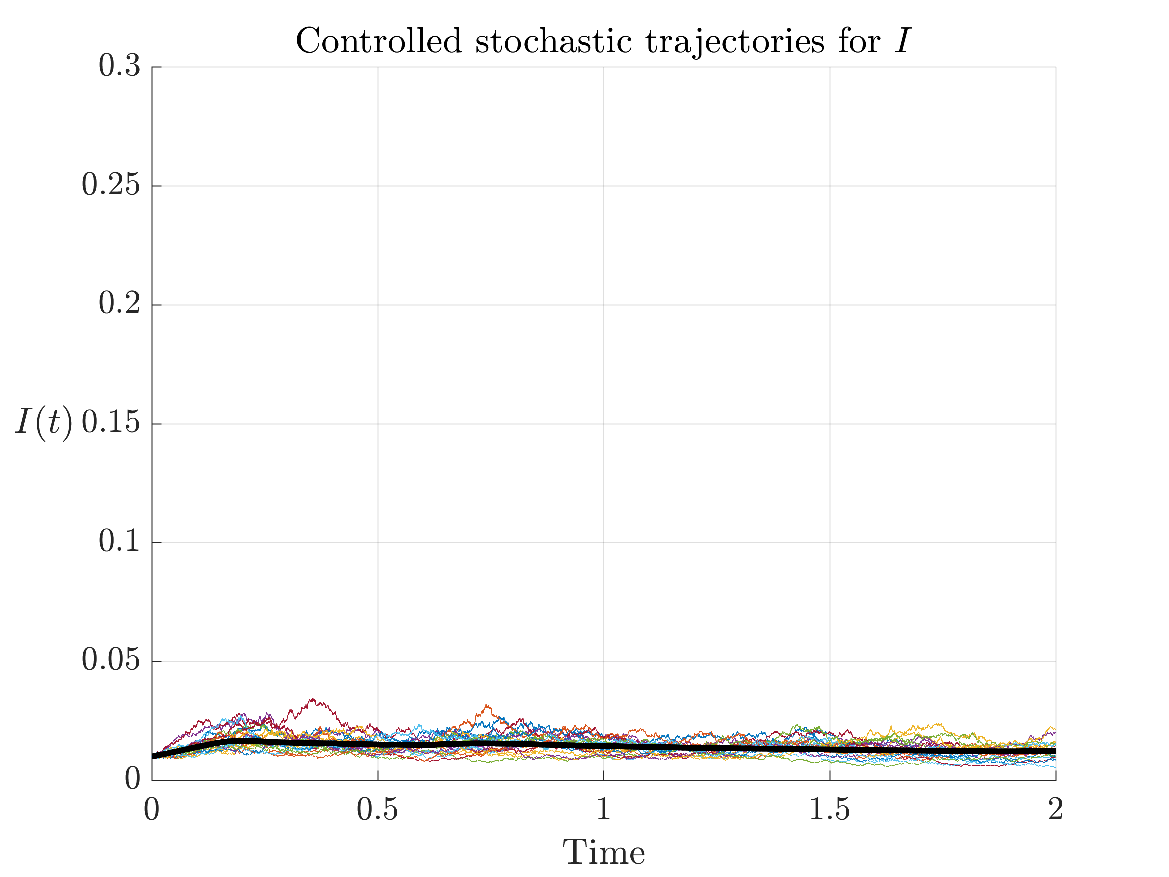}
\includegraphics[width=0.27\textwidth]{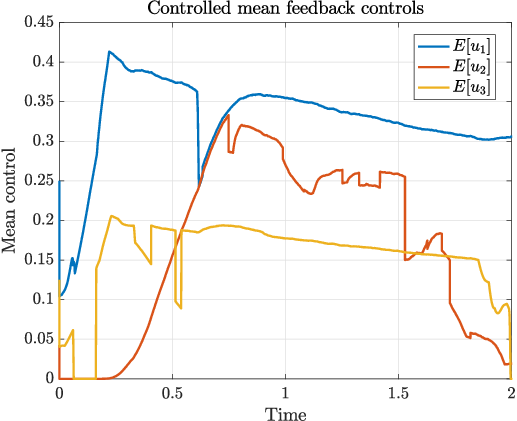}\\[-1mm]

\includegraphics[width=0.31\textwidth]{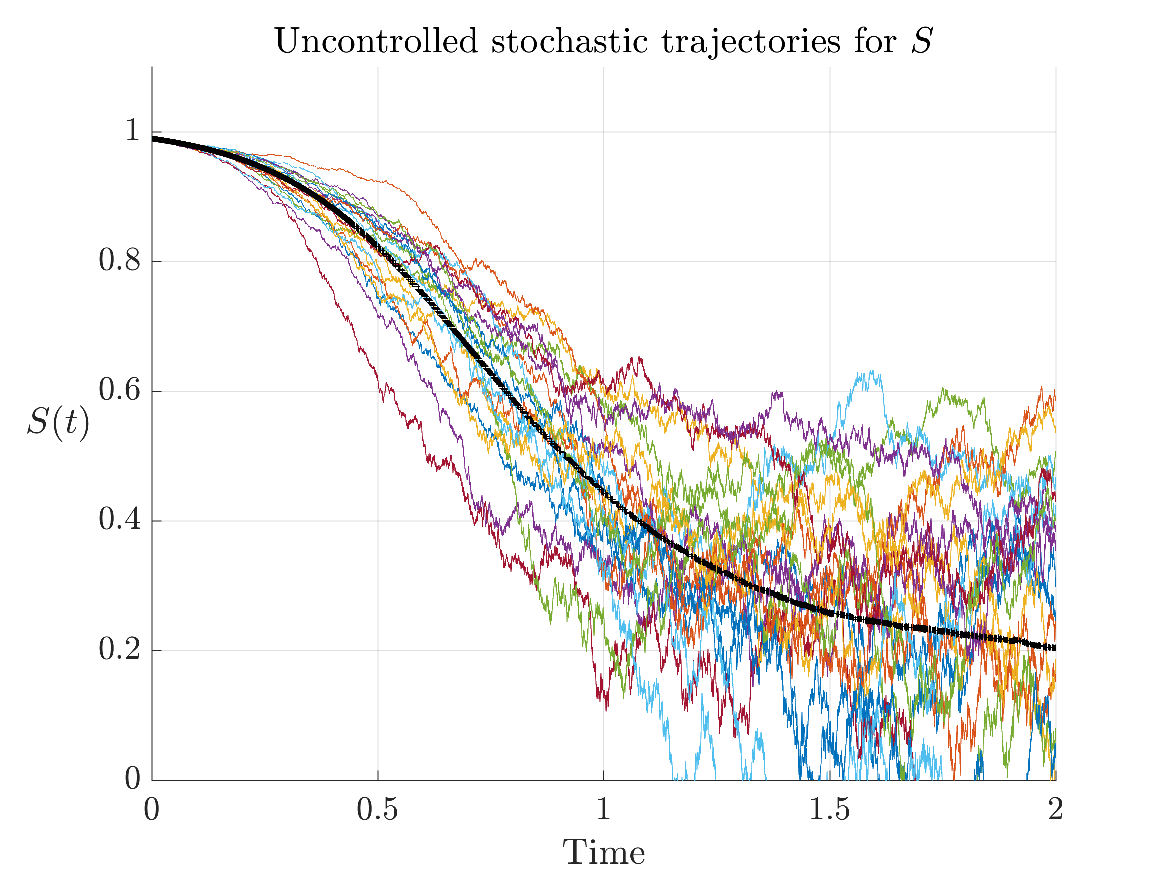}
\includegraphics[width=0.31\textwidth]{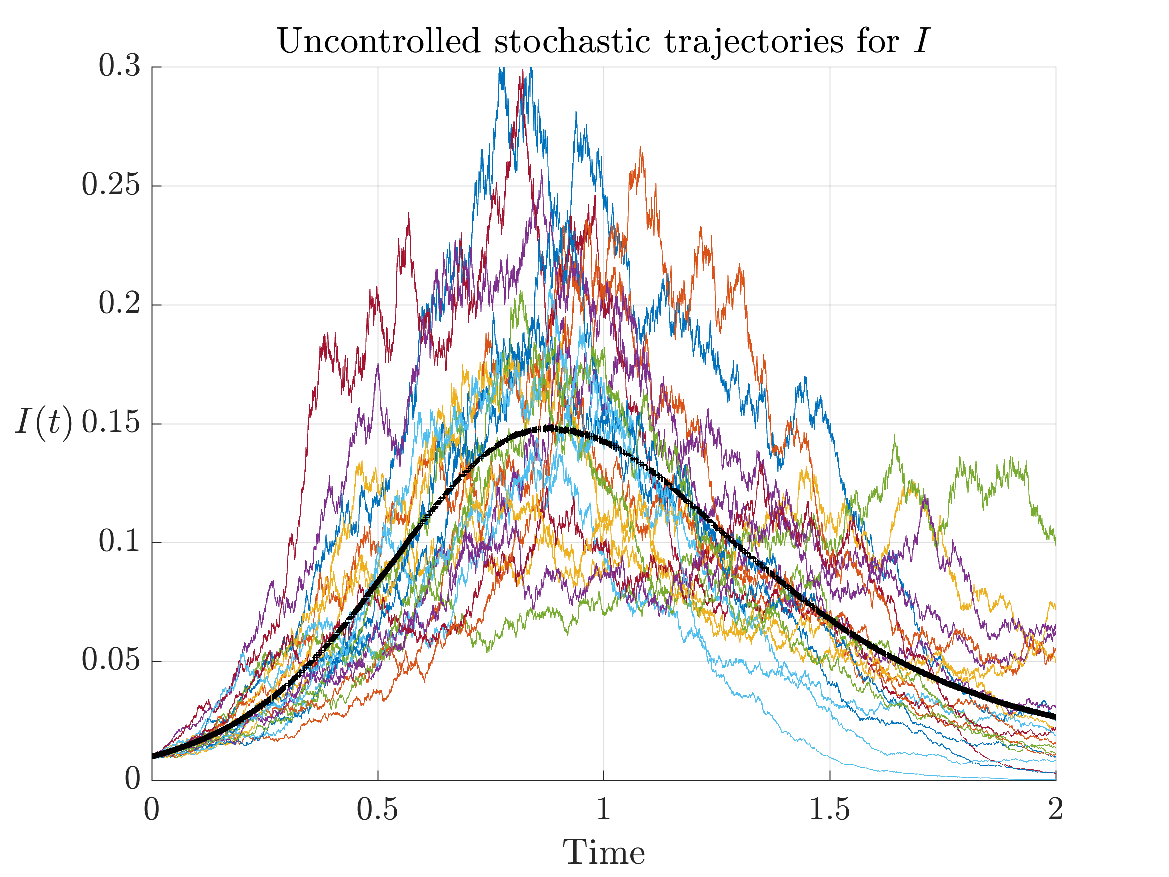}
\includegraphics[width=0.27\textwidth]{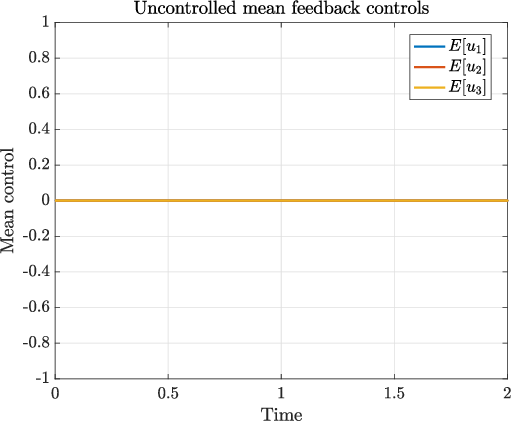}

\caption{Monte Carlo simulations of the controlled and uncontrolled
stochastic system. The top row shows the controlled stochastic
trajectories for the susceptible state $S$ and infected state $I$,
together with the mean feedback controls evaluated along the stochastic
trajectories. The bottom row shows the corresponding uncontrolled
stochastic trajectories and mean controls. The individual curves
represent stochastic realizations, while the emphasized curve denotes
the ensemble mean.}
\label{fig:controlled_uncontrolled_SDE_Case3}
\end{figure}

The corresponding controlled and uncontrolled stochastic simulations
are shown in
Figure~\ref{fig:controlled_uncontrolled_SDE_Case3}.  In the uncontrolled
system, transmission produces a substantial decline in the susceptible
population, whereas the optimized feedback policy limits this loss and
maintains the controlled trajectories at larger values of $S$. This
behavior is consistent with the terminal reward, which favors
realizations that reach the final time with $S>\overline S$. Although the
objective contains no direct running penalty on $I$, the controlled
infected trajectories are also affected because the interventions alter
the transitions among the susceptible, infected, and recovered
populations. In particular, $u_1$ suppresses the $S\to I$ transition,
$u_2$ suppresses the $R\to S$ transition, and $u_3$ enhances the
$I\to R$ transition. The mean controls along the stochastic
trajectories consequently illustrate a policy designed primarily to
preserve the susceptible population at the terminal time, rather than
to minimize the infected fraction throughout the entire time horizon.\\

\begin{figure}[h!]
\centering

\includegraphics[width=0.31\textwidth]{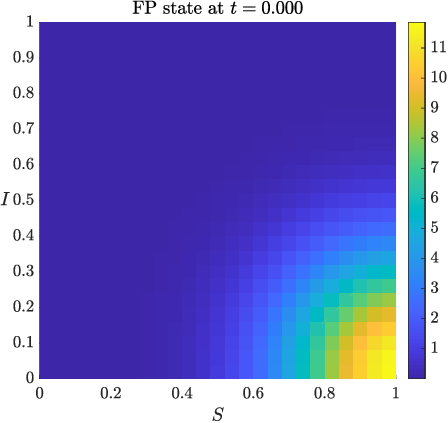}
\includegraphics[width=0.31\textwidth]{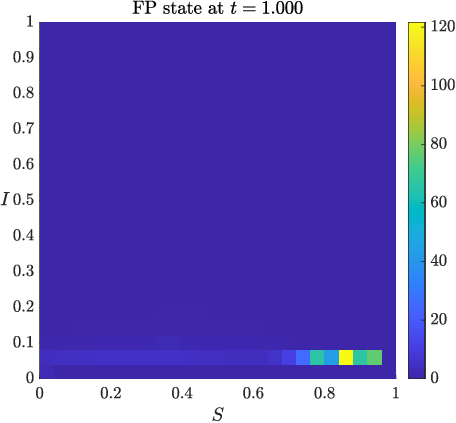}
\includegraphics[width=0.31\textwidth]{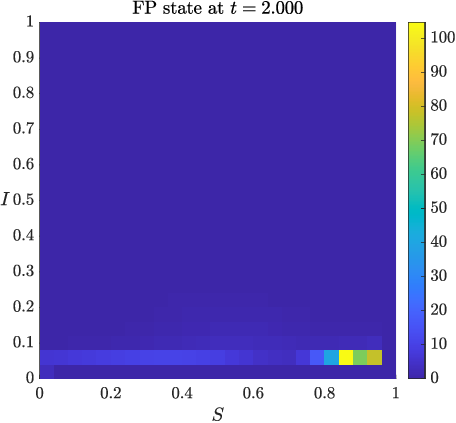}\\

\includegraphics[width=0.31\textwidth]{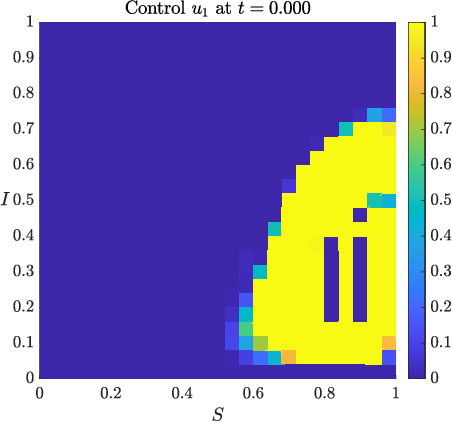}
\includegraphics[width=0.31\textwidth]{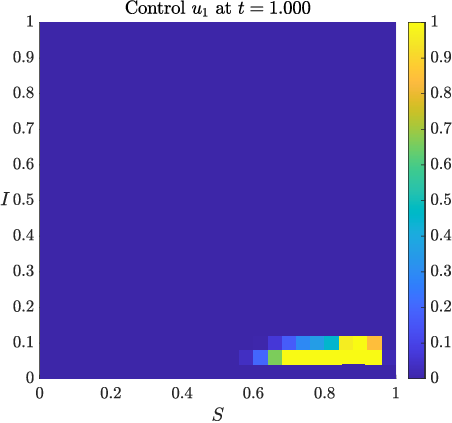}
\includegraphics[width=0.31\textwidth]{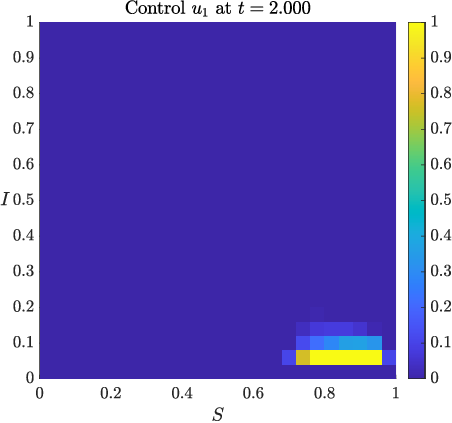}\\

\includegraphics[width=0.31\textwidth]{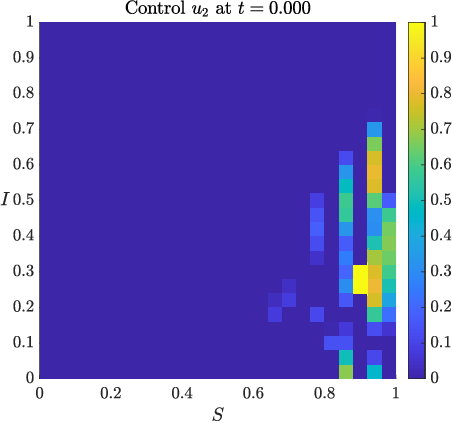}
\includegraphics[width=0.31\textwidth]{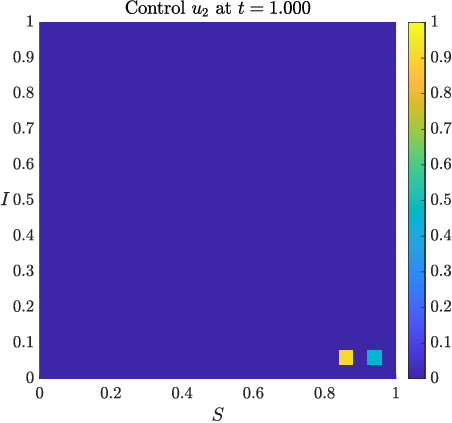}
\includegraphics[width=0.31\textwidth]{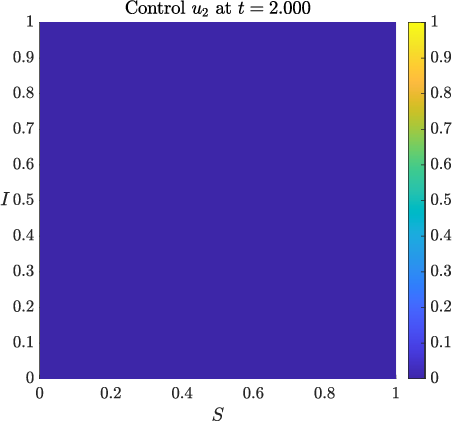}\\

\includegraphics[width=0.31\textwidth]{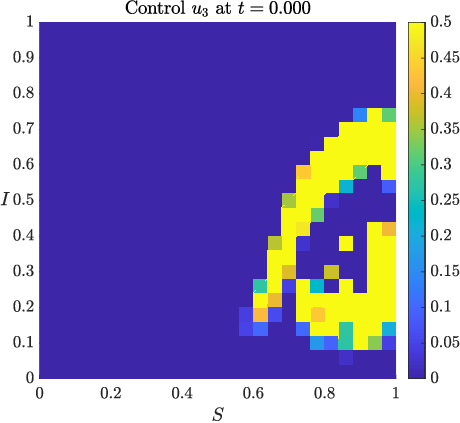}
\includegraphics[width=0.31\textwidth]{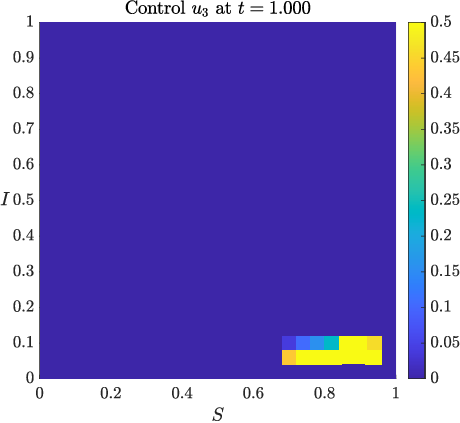}
\includegraphics[width=0.31\textwidth]{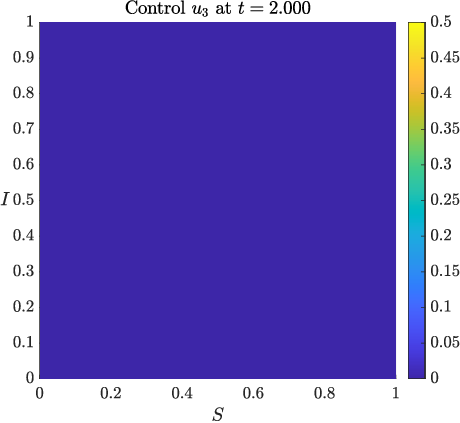}

\caption{Controlled FP solution and optimal feedback controls.
The columns correspond, from left to right, to the initial time $t=0$,
the midpoint $t=T/2$, and the final time $t=T$. The first row shows the
controlled probability density $f(S,I,t)$, while the second, third, and
fourth rows show the optimal feedback controls $u_1(S,I,t)$,
$u_2(S,I,t)$, and $u_3(S,I,t)$, respectively.}
\label{fig:controlled_FP_controls_Case4}
\end{figure}

\noindent {\bf Scenario 4.}
In this scenario, we consider
$
G(S,I,t)=\mathbf{1}_{\{I\geq 0.08\}}(I-0.08),
~ K(S,I)=-\max(S-\overline S,0),
$
combining a threshold-based running penalty on infection with a terminal
reward for maintaining the susceptible fraction above the prescribed
level $\overline S$. The controlled FP density and optimal feedback controls are
shown in Figure~\ref{fig:controlled_FP_controls_Case4}.  In this case,
the feedback policy must balance two state objectives: limiting the
probability of occupying states with $I>0.08$ during the evolution and
steering the terminal density toward states with $S>\overline S$. The control
$u_1$ is particularly relevant to both objectives because reducing
transmission suppresses the $S\to I$ transition, simultaneously limiting
the growth of infection and preserving susceptible individuals. The
control $u_2$ reduces the loss of immunity and therefore suppresses the
$R\to S$ transition, while $u_3$ increases recovery and enhances the
$I\to R$ transition. The different state- and time-dependent structures
of the three feedback maps reflect how the optimizer combines these
mechanisms to balance suppression of infection above the threshold with
the terminal susceptible-population objective.

\begin{figure}[h!]
\centering

\includegraphics[width=0.31\textwidth]{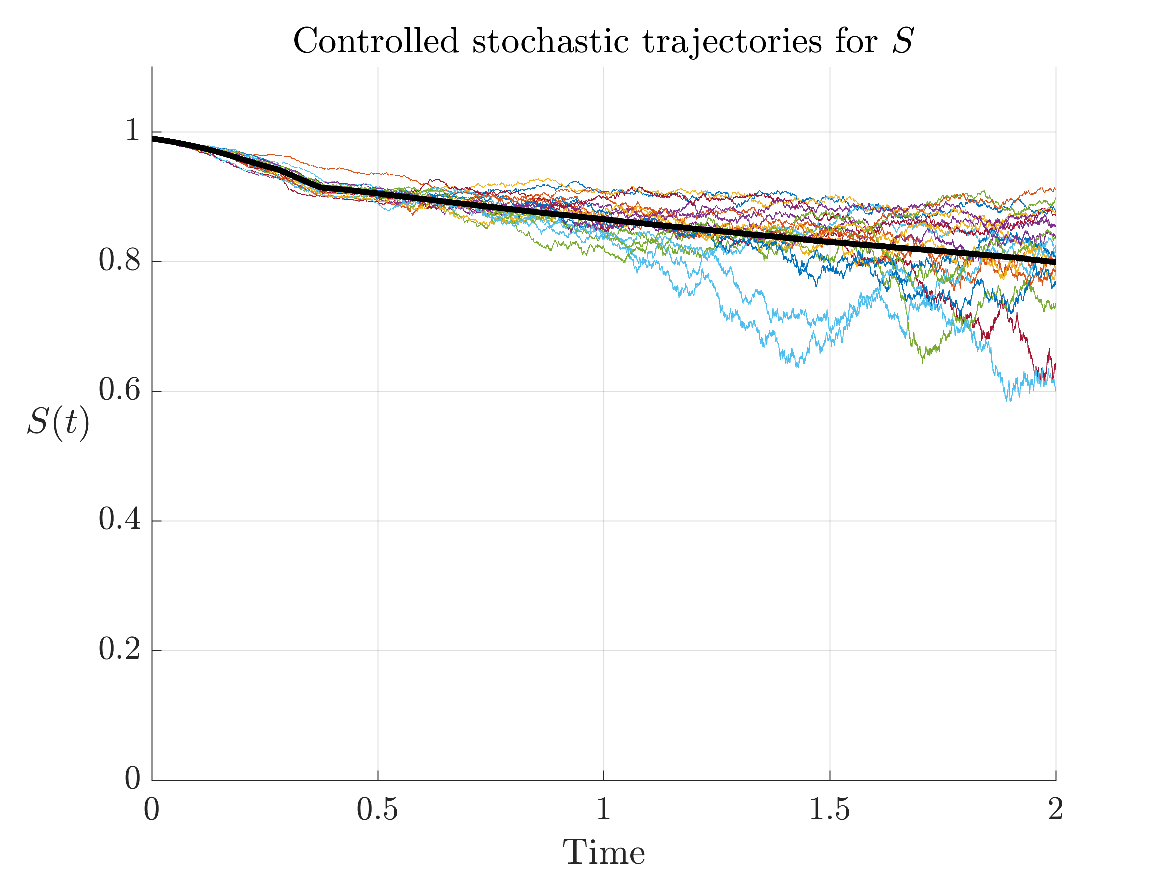}
\includegraphics[width=0.31\textwidth]{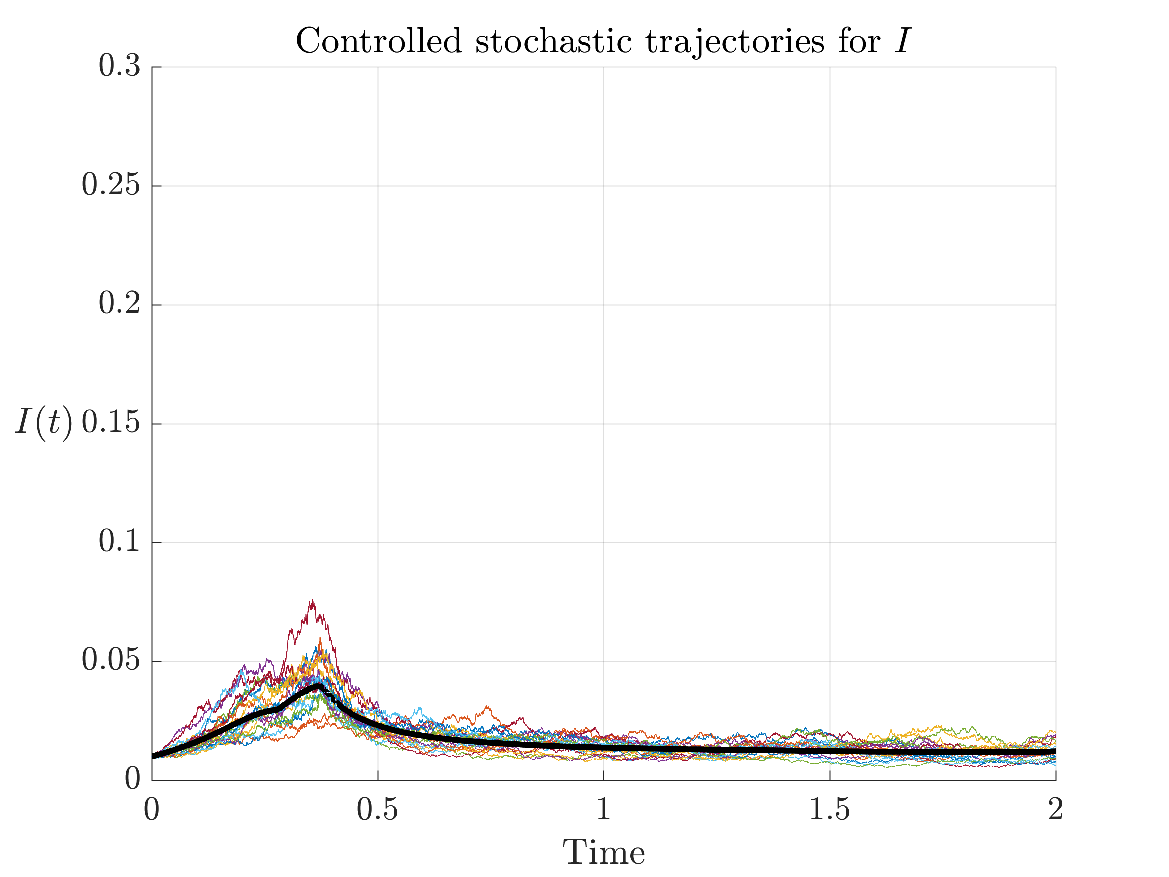}
\includegraphics[width=0.27\textwidth]{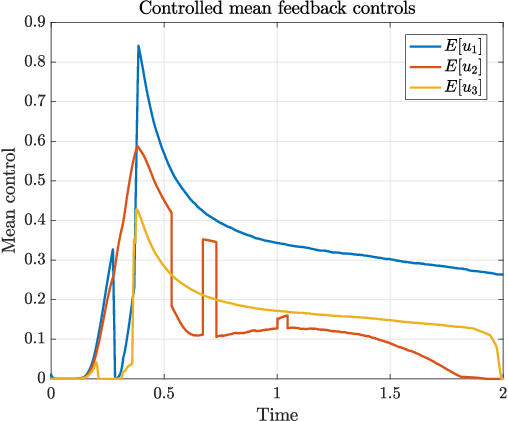}\\[-1mm]

\includegraphics[width=0.31\textwidth]{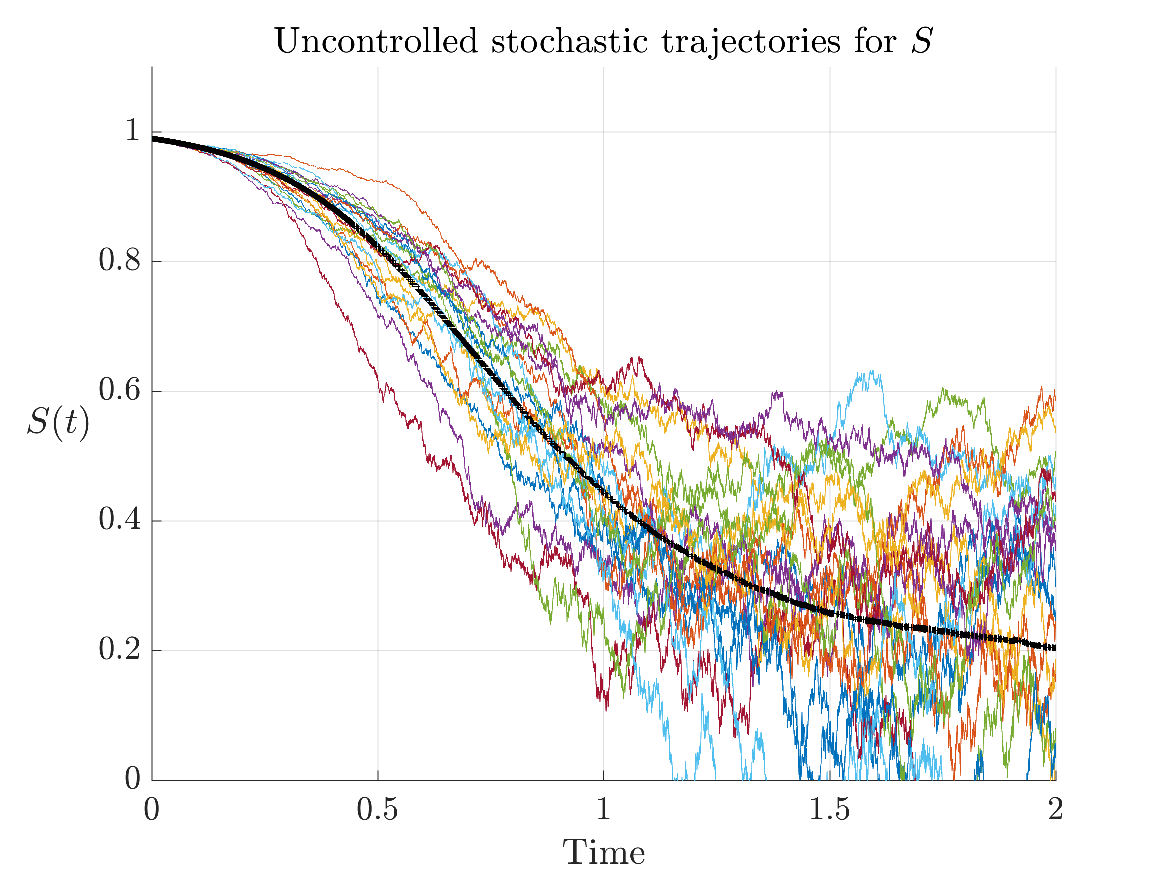}
\includegraphics[width=0.31\textwidth]{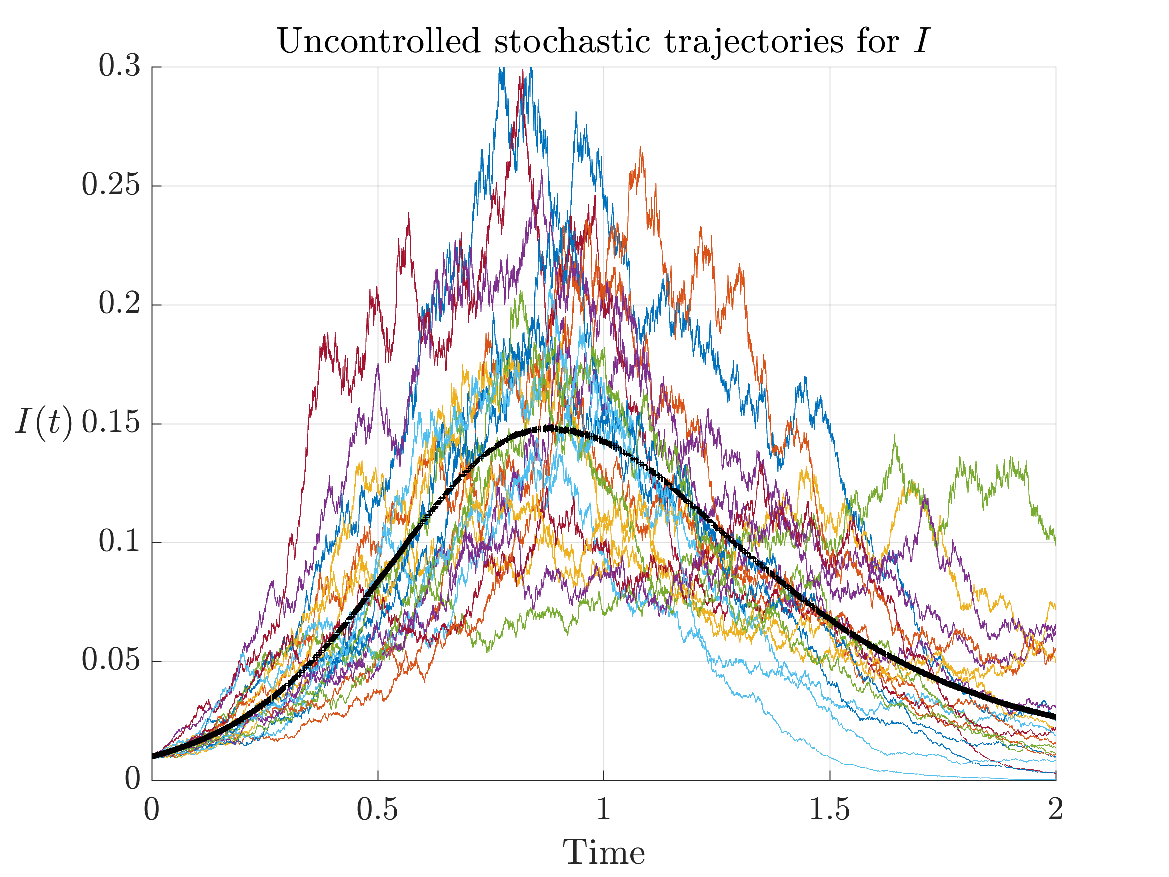}
\includegraphics[width=0.27\textwidth]{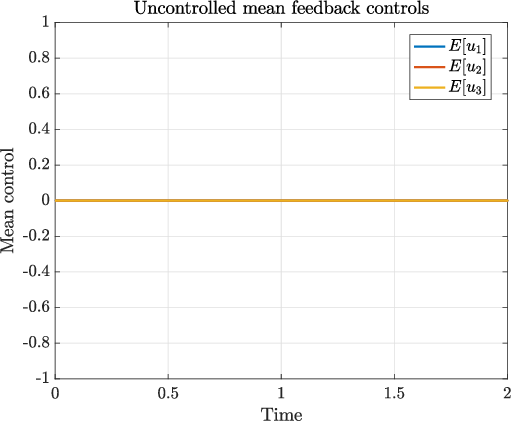}

\caption{Monte Carlo simulations of the controlled and uncontrolled
stochastic system. The top row shows the controlled stochastic
trajectories for the susceptible state $S$ and infected state $I$,
together with the mean feedback controls evaluated along the stochastic
trajectories. The bottom row shows the corresponding uncontrolled
stochastic trajectories and mean controls. The individual curves
represent stochastic realizations, while the emphasized curve denotes
the ensemble mean.}
\label{fig:controlled_uncontrolled_SDE_Case4}
\end{figure}

The corresponding controlled and uncontrolled stochastic simulations
are shown in
Figure~\ref{fig:controlled_uncontrolled_SDE_Case4}. 
 The uncontrolled
trajectories exhibit a substantial decrease in the susceptible fraction
together with an increase in infection that carries the population into
the penalized region $I>0.08$. Under the optimized feedback policy, the
growth of the infected population is reduced while a larger susceptible
fraction is maintained, reflecting the simultaneous influence of the
running and terminal objectives. The mean feedback controls along the
stochastic trajectories show how the three intervention mechanisms
contribute to this behavior: $u_1$ reduces the generation of new
infections and the associated loss of susceptible individuals, $u_2$
retains recovered individuals in the immune class by reducing the
$R\to S$ transition, and $u_3$ accelerates recovery through the
$I\to R$ transition. Thus, Scenario~4 illustrates how the feedback
framework balances short-term regulation of excursions above the
infection threshold with the terminal objective of maintaining a
sufficiently large susceptible population.


\section{Conclusion}\label{sec:conclusion}

In this work, we developed a FP optimal control framework
for a stochastic SIRS
epidemic model with three intervention mechanisms targeting disease
transmission, loss of immunity, and recovery. By incorporating
space- and time-dependent feedback controls into both the drift and
diffusion coefficients, the proposed approach enables interventions
to influence the probability distribution of epidemic outcomes rather
than individual trajectories alone. We investigated the existence of
optimal controls and derived first-order necessary optimality
conditions involving a forward Fokker--Planck equation, a backward
adjoint equation, and a Hamiltonian minimization condition. For
numerical computation, we combined a Chang--Cooper discretization
with a finite-difference treatment of cross-diffusion and a SQH algorithm, together with an analysis of
its convergence properties. Numerical experiments involving four
epidemiological objectives demonstrated how the resulting feedback
policies adapt to different goals, including reducing cumulative
infection, limiting excursions above an infection threshold, and
maintaining a desired susceptible population at the terminal time.
Monte Carlo simulations further illustrated the effectiveness of
the optimized controls relative to uncontrolled stochastic dynamics.
Future work will focus on extending the framework to more complex
epidemic models, strengthening the analysis and numerical treatment
of degenerate diffusion, and incorporating epidemiological data
for model calibration and validation.

\section*{Acknowledgements} S. Roy was supported by the US National Science Foundation grant nos. 2212938, 2309491, 2230790.

\bibliographystyle{plain}
\bibliography{ref}

\end{document}